\documentclass[11pt,a4paper,reqno,final]{article}
\usepackage{template}
\usepackage{xfrac}

\DeclareMathOperator\E{E}
\let\P\relax
\DeclareMathOperator\P{P}

\newcommand\NN{\mathbb{N}}
\newcommand\ZZ{\mathbb{Z}}

\newcommand\RR{\mathbb{R}}

\newtheorem{application}{Application}
\usepackage{stmaryrd}

\usepackage[backend=bibtex,giveninits=true,maxcitenames=12,maxbibnames=9]{biblatex}
\usetikzlibrary{shapes.misc}
\tikzset{cross/.style={cross out, draw=black, minimum size=5*(#1-\pgflinewidth), inner sep=0pt, outer sep=0pt},
cross/.default={1pt}}

\makeatletter
\newcommand\xleftrightarrow[2][]{%
  \ext@arrow 9999{\longleftrightarrowfill@}{#1}{#2}}
\newcommand\longleftrightarrowfill@{%
  \arrowfill@\leftarrow\relbar\rightarrow}
\makeatother

  \newcounter{iconst}

\begin{document}

\title{Random Markov property for random walks in random environments}

\date{\today}
\author{
  Julien Allasia
  \thanks{Email: \texttt{allasia@math.univ-lyon1.fr}; Universite Claude Bernard Lyon 1, CNRS, Ecole Centrale de Lyon, INSA Lyon, Université Jean Monnet, ICJ UMR5208,
69622 Villeurbanne, France.}
  \and
Rangel Baldasso
  \thanks{Email: \ \texttt{rangel@puc-rio.br}; \ Department of Mathematics, PUC-Rio, Rua Marqu\^es de S\~ao Vicente 225, G\'avea, 22451-900 Rio de Janeiro, RJ - Brazil.}
  \and
  Oriane Blondel
  \thanks{Email:  \ \texttt{blondel@lpsm.paris}; \  Universit\'e Paris Cit\'e and Sorbonne Universit\'e, CNRS, Laboratoire de Probabilit\'es, Statistique et Mod\'elisation, F-75013 Paris, France.}
  \and
  Augusto Teixeira
  \thanks{Email: \ \texttt{augusto@impa.br}; \ IMPA, Estrada Dona Castorina 110, 22460-320 Rio de Janeiro, RJ - Brazil and IST, University of Lisbon, Portugal.}
}

\maketitle

\begin{abstract}
We consider random walks in dynamic random environments and propose a criterion which, if satisfied, allows to decompose the random walk trajectory into i.i.d.\ increments, and ultimately to prove limit theorems. The criterion involves the construction of a random field built from the environment, that has to satisfy a certain random Markov property along with some mixing estimates. We apply this criterion to correlated environments such as Boolean percolation and renewal chains featuring polynomial decay of correlations.

  \medskip


\end{abstract}

\section{Introduction}\label{s:intro}
Random walks in random environments have been a subject of intensive study since the 60s-70s, when the first models motivated by applications in biophysics or chemistry were introduced \cite{chernov1967replication,Temkin1972}. One of the main challenges in their study is that memory of the past trajectory of the random walk can be conserved through the environment currently explored. A large part of the literature on the subject focuses on the identification of conditions (on the random walk jumps and the law of the environment) that imply memory loss in some sense. 

One of the most prominent example in this direction is \cite{Sznitman1999}. The environment is assumed to be i.i.d., and the random walk to satisfy some ballisticity condition. This allows to build regeneration times for the trajectory of the random walk { and prove that it satisfies a law of large numbers}. Let us note that it is natural to request some kind of directional transience, since otherwise the random walk keeps revisiting previously explored environment and there is no reason why it should forget what it saw (indeed in general it does not, see \cite{Solomon1975}). The i.i.d.\ assumption on the other hand is clearly non-optimal, and later works have focused on relaxing it. A landmark paper in that direction is \cite{comets2004}, where the authors introduce a \emph{cone-mixing} condition that the environment has to satisfy. It allows to build a sequence of approximate regeneration times (which yield increments that are close to being i.i.d.) { and once again obtains law of large numbers}.

This has later been refined in various ways \cite{Bethuelsen2016, bethuelsen2023limitlawsrandomwalks}, especially in the context of random walks in dynamic random environments (RWDRE). In this case, one may consider the time direction as the direction of transience of a random walk in dimension $d+1$ whose last coordinate deterministically increases by $1$ at each step. Ballisticity is then automatic, and the setting allows to consider weaker decoupling properties on the environment. The main drawback of the cone-mixing property is that is requires a \emph{uniform} decoupling property, which many natural dynamic environments do not satisfy. 

In recent years, there has been a growing interest in studying RWDRE when the environment has a non-uniform decoupling property \cite{Avena2017}, especially using renormalization techniques \cite{cavalinhos, HKT19, baldasso2023fluctuationboundssymmetricrandom, arcanjo2023lawlargenumbersballistic, allasia2023law, allasia2024asymptoticdirectionballisticrandom, berard2016fluctuations,den2014random}. { These works examine different aspects of the behavior of the random walk under distinct dynamic environments, with main focus given to the establishment of laws of large numbers and central limit theorems.} What we propose in the present paper is to give a framework in which one can build a sequence of regeneration times for the RWDRE, and use renormalization to control the tails of these regeneration times in terms of decorrelation properties of the environment. We can then apply standard limit results for sums of i.i.d.\ random variables to the trajectory of the RWDRE.

It is worth mentioning that the same type of problems have been studied by means of analyzing the environment seen from the particle~\cite{dolgopyat2008random, dolgopyat2009non, redig2013random,Avena2017,bethuelsen2023limitlawsrandomwalks}.

Our strategy relies on the construction of a field $\eta\in\{0,1\}^{\Z^{d+1}}$, which is model-dependent (see examples in Section~\ref{s:examples}). $\eta$ has to be built from the environment (and possibly extra independent randomness) in such a way that it has the \emph{Random Markov Property} which we define below (Definition~\ref{d:rmp}). Intuitively, the Random Markov Property says that, conditional on $\eta_x=1$, the future of a random walk trajectory going through $x$ is independent of its past; this will allow us to decompose the trajectory of the RWDRE into i.i.d.\ increments (Lemma~\ref{l:indep_decomp}), conditioning on the event that $\eta$ is $1$ at the origin. Note that this step is less straightforward than it might appear from the intuitive interpretation of the Random Markov Property. For instance, it necessitates resampling the field $\eta$ (Lemma~\ref{l:resampling_lemma}) in order to erase irrelevant information. The tail of the regeneration times is controlled by the decorrelation properties of $\eta$ (see Definition~\ref{d:decoupling} and Lemma~\ref{l:moments}), which in turn are inherited from those of the environment in our examples. The tail estimates are proved through a renormalization procedure.

The paper is organised as follows. In Section~\ref{ss:rwdre}, we give a short definition of the random walk in dynamic random environment. In Section~\ref{ss:rmp}, we explain what we need from the auxiliary field $\eta$, and in particular the Random Markov Property (RMP). We conclude Section~\ref{s:intro} with the statement of our main result Theorem~\ref{t:main_theorem}. In Section~\ref{s:graphical} we give a formal construction of our random walk and the surrounding objects. Section~\ref{s:ind_dec} is devoted to the construction of regeneration times for the trajectory of the random walk, using the RMP and a resampling procedure. Section~\ref{s:renewal} presents the renormalization strategy that allows to control the tails of the regeneration times, and the proof of our main results is given in Section~\ref{ss:final_proof}. We conclude in Section~\ref{s:examples} with examples of dynamic random environments for which we construct an appropriate field $\eta$: Boolean percolation with polynomial tails for the radii, and independent renewal chains with polynomial tails for the interarrival distribution.

\subsection{Notation}

For $x\in\Z^d$ and $t\in\Z$, we see $z=(x,t)$ as an element of $\Z^{d+1}.$ We denote by $|\cdot|$ the $L^\infty$ norm on $\Z^{d+1}$.
Mind that depending on the context, $0$ can mean the origin of $\mathbb{Z}^d$ or $\mathbb{Z}^{d + 1}$. For instance, when we write $(0, t)$ with $t\in\ZZ$ and a point in $\mathbb{Z}^{d + 1}$ was expected, what we mean is $0 \in \mathbb{Z}^d$.

$\NN$ denotes the set of natural integers starting from $0$. $\NN^*$ is $\NN\setminus\{0\}$.
For any two integers $a < b$, denote by $\llbracket a,b\rrbracket^d = [a,b]^{d} \cap \Z^{d}$.

\subsection{Random walk in dynamic random environment}\label{ss:rwdre}

We start by defining somewhat informally the environment and random walk that we consider. In Section~\ref{s:graphical} we propose a formal construction of those objects.

Fix $S$ a countable set and $R$ a positive integer. We work with an environment $\omega=(\omega_{z})_{z\in\Z^{d+1}}\in S^{\Z^{d+1}}$ that will be chosen with some fixed probability distribution $\P$. We assume throughout the text that $S$ is endowed with a distance that makes it a Polish metric space.

To define the random walk in random environment, for each state $s \in S$ we choose a probability kernel $p(s, \cdot)$ on $\llbracket -R,R\rrbracket^d$. For $\omega\in \Omega$, the random walk in the environment $\omega$ is the process $X=(X_t)_{t\in\N}$ such that $X_0=0$ and, for $t\in\NN$, $x\in\ZZ^d$ and $y\in\llbracket -R,R\rrbracket^d$,
\begin{equation}\label{e:jump_law}
\P^\omega(X_{t+1}=x+y | X_t=x) = p(\omega_{(x,t)}, y).
\end{equation}
We say that $R$ is the \emph{range} of $X$, and that $X$ is \emph{finite-range}. 

\nc{c:ellipticity}
We assume that $X$ is \emph{uniformly elliptic}, meaning that
\begin{equation}
  \label{e:ellipticity}
  \min_{y \in \llbracket -1, 1 \rrbracket^d} \inf_{s \in S} p(s,y) = \uc{c:ellipticity} > 0.
\end{equation}
Mind that although $X$ has range $R$, we only ask for uniform {ellipticity} for the jumps in $\llbracket -1,1\rrbracket^d$.

We see the trajectory of our random walk as a path in $\Z^{d+1}$ by considering $Z=(Z_t)_{t\in\N}$, where
\begin{equation}
 Z_t=(X_t,t).
\end{equation}

\subsection{The field $\eta$: Random Markov Property and decoupling}\label{ss:rmp}

Our main assumption is the existence of a random field $\eta=(\eta_z)_{z\in\Z^{d+1}}:\Z^{d+1}\rightarrow \{0,1\}$ on the same probability space as $\omega$, which will help us construct a renewal structure based on times where the random walk hits a point where $\eta_z=1$ (note that we will keep the notation $\P$ for the probability space of the environment, but now including not only $\omega$ as before, but also $\eta$). Let us start by some technical assumptions.

  \nc{c:eta_zero}
\begin{assumption}\label{a:assumptions_eta} We ask $\omega$ and $\eta$ to satisfy the following conditions.
\begin{enumerate}
 \item $(\omega,\eta)$ is translation-invariant: for any $z\in\ZZ^{d+1}$, the law of $(\omega_{z+\cdot},\eta_{z+\cdot})$ under $\P$ is the same as that of $(\omega,\eta)$.
\item There exists $\uc{c:eta_zero}>0$ such that
\begin{equation}
    \label{e:eta_zero}
    \P(\eta_{0}=1)\geqslant \uc{c:eta_zero}.
  \end{equation}
\item There exists $s_\star\in S$ such that
    \begin{align}\label{e:condition_eta}
     \forall z\in\ZZ^{d+1},\;\eta_z=1\Rightarrow \omega_z=s_\star.
    \end{align}
\end{enumerate}
\end{assumption}

The translation invariance is a natural requirement in our setting. The second assumption will be used to ensure that the random walk often meets points where $\eta$ is $1$. The third assumption is a technicality that helps us remove an extra condition in the construction of our renewal times. { Without this last assumption, the field $\eta$ would provide a renewal of the environment, but not necessarily of the random walk, as it still needs to observe the environment at its current position, that might depend on the cone of the past.}
Bear in mind that in our examples, $\eta$ will be a \emph{non-local} function of $\omega$ (and sometimes requiring extra randomness).

{\begin{remark}
 In the following we may want to incorporate extra randomness into $\omega$ in order to allow us to build $\eta$. Note however that the probability kernels $p(s,\cdot)$ will not look at that extra information. In particular, we only need \eqref{e:condition_eta} in the sense that the environment as perceived by the random walker takes a specific value $s_\star$.
\end{remark}}

The crucial assumptions on $\eta$ are given by two definitions that we present now. 
For $z=(x,t)\in\Z^{d+1}$, let us define the future and past cones rooted at $z$.

\begin{figure}[h]
    \centering
    \begin{tikzpicture}[use Hobby shortcut,scale=0.7]
        \fill[black!10] (-4,-2)--(0,0)--(4,-2);
        \fill[black!10] (-4,2)--(0,0)--(4,2);
        \draw[right] (0,-1) node {$C_0^-$};
        \draw[left] (0,1) node {$C_0^+$};
        \draw[->,>=latex, right] (-4,0)--(4,0) node {$x$};
        \draw[->,>=latex,above] (0,-2.2)--(0,2.4) node {$t$};
        \draw (-4,-2)--(4,2);
        \draw (-4,2)--(4,-2);
    \end{tikzpicture}
    \caption{Future and past cones rooted at the origin when $d=1$ and $R=2$.}
    \label{f:decoupling}
\end{figure}
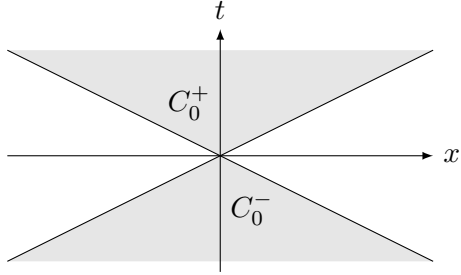

\begin{equation}\begin{split}
C_0^{+}&=\{(x,t)\in\Z^{d+1} : |x|\leqslant Rt\},\\
C_0^-&=\{(x,t)\in\Z^{d+1} : t\leqslant 0 \text{ and } |x|\leqslant R|t|\},\\
C_z^\pm&=z+C_0^\pm.
                \end{split}
\end{equation}
Note that, if $X_t=x$, since $X$ has range $R$, we have
\begin{equation}
  \begin{split}
    Z_{t+s}\in C_{(x,t)}^+&\text{ for all } s \in \N,\\
    Z_{t-s}\in C_{(x,t)}^-&\text{ for all } s\in\llbracket 0,t\rrbracket.
  \end{split}
\end{equation}
Let us denote, for $z\in\Z^{d+1}$,
\begin{equation}\label{eq:sigma_algebras}
 \begin{split}
   \cF_z^-&=\sigma\left((\omega_{z'},\eta_{z'}): z' \in C^-_z\right),\\
   \cF_z^+&=\sigma\left(\omega_{z'}: z' \in C^+_z\right).
 \end{split}
\end{equation}
Mind the subtlety in the definition of $\cF_z^-$ and $\cF_z^+$: for the former, we look at both processes $\omega$ and $\eta$ in the past cone; for the latter, we only look at $\omega$ in the future cone.

\begin{definition}[Random Markov property]\label{d:rmp}
  We say the field $(\omega,\eta)$ satisfies the Random Markov property if, for any $z\in\Z^{d+1}$ and for any bounded random variable $W$ measurable with respect to $\cF^+_z$,
  \begin{equation}\label{e:rmp}
    \E[W|\cF^-_z]\,\eta_z=\E[W|\eta_z]\,\eta_z.
  \end{equation}
\end{definition}

  \begin{remark}
    It is important to emphasize that $\mathcal{F}^-$ depends on both $\omega$ and $\eta$ in the past cone $C^-$, while $\mathcal{F}^+$ depends only on $\omega$ in $C^+$.
    This is a crucial difference that makes the proof more complicated in many places.
    In other words, it would have been much easier to prove a CLT under the hypothesis \eqref{e:rmp} if $\mathcal{F}^+$ was allowed to depend on $\eta$ as well.
    However, this would rule out several possible applications of the model, as in the case of independent renewal chains that we present below.
  \end{remark}
  

Additionally, we need a decoupling property for the field $\eta$. 
Let $A\subseteq\Z^{d+1}$. We say that $A$ is a box if it is of the form $\prod_{i=1}^{d+1}\llbracket a_i,b_i\rrbracket$ with $a_i\leqslant b_i$ integers. If $A$ is a box, we define its spatial diameter $r(A)$ and height $h(A)$ by
\begin{equation}
  \label{eq:diam_h}
  \begin{split}
    r(A)&:=\max_{i=1,\ldots,d}(b_i-a_i),\\
    h(A)&:=b_{d+1}-a_{d+1}.
  \end{split}
\end{equation}
For $A=\prod_{i=1}^{d+1}\llbracket a_i,b_i\rrbracket,A'=\prod_{i=1}^{d+1}\llbracket a'_i,b'_i\rrbracket$ two boxes, we call $\sep(A,A')$ their vertical separation defined by
\begin{equation}
  \label{eq:separation}
  \sep(A,A') =
  \begin{cases}
    a_{d+1}-b'_{d+1}, & \text{ if } b'_{d+1}< a_{d+1},\\
    a'_{d+1}-b_{d+1}, &\text{ if }b_{d+1}<a'_{d+1},\\
    0, & \text{ else.}
  \end{cases}
\end{equation}

\nc{c:decoupling}

\begin{definition}[Decoupling property]\label{d:decoupling}
 Let $\alpha\geqslant 0$. We say that $\eta$ satisfies the $\alpha$--decoupling property if there exists $\uc{c:decoupling}=\uc{c:decoupling}(\alpha)>0$ such that for any $r\in\N$, for any boxes $A,B\subseteq\ZZ^{d+1}$ satisfying 
 \begin{align*}
r(A)\vee r(B)\leqslant (2R+1)r,\quad h(A)\vee h(B)\leqslant r,\quad \sep(A,B)\geqslant r,
\end{align*}
and for any $f_1:\{0,1\}^A\rightarrow\{0,1\}$, $f_2:\{0,1\}^B\rightarrow\{0,1\}$, we have
\begin{align*}
\E[f_1(\eta_A)f_2(\eta_B)]\leqslant \E[f_1(\eta_A)]\E[f_2(\eta_B)]+\uc{c:decoupling} r^{-\alpha}.
\end{align*}
\end{definition}

\begin{remark}\label{r:decoupling_function}
In most cases, such as those presented in Section~\ref{s:examples}, $\eta$ will satisfy a stronger decoupling property that works for larger sets of boxes. We use the following vocabulary in that context. Let $\eps:\R_+^3\rightarrow\R_+$. We say that $\eps$ is a decoupling function for $\eta$ if for any $r,h,s\in\N$, for any boxes $A,B\subseteq\Z^{d+1}$ such that
\begin{align}\label{e:distance_decoupling}
  r(A)\vee r(B)\leqslant r,\quad h(A)\vee h(B)\leqslant h,\quad \sep(A,B)\geqslant s,
\end{align}
and for any $f_1:\{0,1\}^A\rightarrow\{0,1\}$, $f_2:\{0,1\}^B\rightarrow\{0,1\}$, we have
\begin{align*}
  \E[f_1(\eta_A)f_2(\eta_B)]\leqslant \E[f_1(\eta_A)]\E[f_2(\eta_B)]+\eps(r,h,s).
\end{align*}
Observe that, if there exists $c_2 > 0$ such that $\eps((2R+1)r,r,r)\leq c_2r^{-\alpha}$ for any $r\in\N$, then the environment satisfies the $\alpha$--decoupling property.

\end{remark}

\subsection{Main results}

Our main theorem is a law of large numbers and an annealed central limit theorem, meaning that it holds for the annealed law defined by $\bP=\int \P^\omega(\cdot)\,\mathrm{d}\P(\omega)$ (see Section 2).

\begin{theorem}\label{t:main_theorem}
  Assume that $X$ is finite-range and uniformly elliptic and that Assumption \ref{a:assumptions_eta} as well as the Random Markov property and the $\alpha$--decoupling property is satisfied.
  \begin{enumerate}
  \item Law of large numbers. If $\alpha>1$, then, as $t$ goes to infinity, $\frac{X_t}{t}$ converges $\bP$-almost surely to a certain limit $v\in\RR^{d}$.
  \item Central limit theorem. If $\alpha>2$, then, as $t$ goes to infinity, $\frac{X_t - vt}{\sqrt{t}}$ converges in $\bP$-distribution to a Gaussian random variable in $\R^d$ with zero mean value.
  \end{enumerate}
\end{theorem}

The proof of Theorem \ref{t:main_theorem} is divided into two parts, following the road map of the renewal method from \cite{Sznitman1999}. The first one, namely Section \ref{s:ind_dec}, consists in defining a sequence of renewal times $(T_k)_{k\in\NN}$ and show that conditioned on $T_0=0$, the increments of the random walk along this sequence of times {are} i.i.d. This is a consequence of the Random Markov property. In the second part of the proof, namely Section \ref{s:renewal}, we need to control the moments of $T_1$ conditioned on $T_0=0$. Here, the mixing behavior that we demand for $\eta$ will be instrumental. The final proof of Theorem \ref{t:main_theorem} can be found in Section \ref{ss:final_proof}. A construction will be provided for {$v$} and the variance of the limiting Gaussian distribution. In Section \ref{s:examples}, we give examples of environments that satisfy our assumptions.

{\begin{remark}\label{rem:covariance}
 The value $a$ and the variance $ \Sigma $ of the limiting Gaussian distribution are defined in terms of the increment of the random walk between two renewal times. In particular, a simple condition that is sufficient to ensure that $\Sigma$ is non-degenerate is the following.
 \begin{equation}
  \hat \P(\eta_{(e_i,1)}=1)>0 \quad \forall i=1,\ldots,d.
 \end{equation}
Indeed, uniform ellipticity and our definition of renewal times \eqref{e:Y_k} then ensure that with positive probability the increment can take values in all space directions. This condition is in particular satisfied in the examples of Section~\ref{s:examples}.
\end{remark}}

{
\paragraph{Relations to previous work}

Although we already gave a brief overview of some works that have been published on this topic, let us now explore this discussion further, addressing how our article interacts with earlier studies.
One way to describe this relation is through the lens of proof techniques.

The first ingredient in our proof is the construction of renewal times for the random walk.
This technique appeared earlier in many works \cite{kuczek1989central,MR2198849,zbMATH06514478}.
But the main distinction to our case is the fact that one would previously find these renewal times dynamically as one follows/builds the trajectory of the random walk.
In our paper, they are instead constructed at the space-time environment alone (through the field $\eta$) and we later prove that the random walk will hit these points often.
This is the ingredient that makes it easier to find renewals across different environments.
However, our work does not supersede the existing articles that employ renewal times, since there are many model-specific choices that each work employs during these constructions.

Another technique that we employ is renormalization, more precisely to prove that the random walk is very likely to hit the traps given by $\eta$.
Again, this has been applied many times in the past and some of the closest examples to what we do here are \cite{cavalinhos,HKT19}.
But again, our no work does not contain the others, as the techniques and hypotheses vary quite a lot within each paper.

In short, even though we provide a somewhat general technique for establishing CLT, applying our theorem for each example of environment requires a specific understanding of the model in question and may not be easy to do, as we discuss in the Open Problems' session.
}

\paragraph{Acknowledgements.}

The research of JA was supported by a doctoral contract and travel funds provided by CNRS Mathématiques.
RB has counted on the support of ``Conselho Nacional de Desenvolvimento Científico e Tecnológico – CNPq'' grants ``Produtividade em Pesquisa'' (308018/2022-2) and ``Projeto Universal'' (402952/2023-5).
The research of OB was partially supported by the ANR grant MICMOV (ANR-
19-CE40-0012) of the French National Research Agency (ANR).
During this period, AT has also been supported by grants ``Projeto Universal'' (406250/2016-2) and ``Produtividade em Pesquisa'' (304437/2018-2) from CNPq and ``Jovem Cientista do Nosso Estado'' (202.716/2018) from FAPERJ.

\section{Formal construction}
\label{s:graphical}

Let us endow $S^{\ZZ^{d+1}}$ with the product topology and the subsequent Borel $\sigma$-algebra. Let $(\Omega,\mathcal{T},\P)$ be a probability space and $\omega,\eta$ two random variables on $(\Omega,\mathcal{T})$ taking values in $S^{\ZZ^{d+1}}$ and $\{0,1\}^{\ZZ^{d+1}}$ respectively. { The random field $\eta$ is allowed to depend both on $\omega$ and on additional randomness. By possibly incorporating this randomness into $\omega$, we will also denote elements of $\Omega$ by $\omega$.} We denote by $\hat{\P}$ the probability distribution on $(\Omega, \mathcal{T})$ given by $\P$ conditioned on $\eta_{0}=1$ (recall \eqref{e:eta_zero}), that is, $\hat{\P}(A) = \P (A | \eta_{0}=1)$.

We now assume that environment $\omega$ is fixed and provide a formal construction of the random walk in $\omega$. What we want to do, besides having equality \eqref{e:jump_law} satisfied, is to couple random walks starting from different points together. In order to do so, we use uniform random variables to encapsulate the randomness needed for the jumps of the random walks and allocate them to points in $\ZZ^{d+1}$ instead of allocating them to the random walks themselves, depending on their starting points.

For each $s \in S$, partition $[0,1]$ into $(2R+1)^d$ intervals $(I_s^y)_{y\in \llbracket -R,R\rrbracket^d}$ so that the length of $I_s^y$ is $p(s,y)$. Set, for $u\in [0,1]$,
\begin{equation}
g(s,u)=\sum_{y\in\llbracket -R,R\rrbracket^d} y \mathbf{1}_{I_s^y}(u).
\end{equation}
This is a measurable function from $S \times [0,1]$ to $\llbracket -R,R \rrbracket^d$ that will determine the jump of a random walk if the state of the environment at its location is $s$.

In what follows we will augment the probability spaces and distributions $\P^\omega$, $\mathbb{P}$ and $\hat{\mathbb{P}}$ (which previously denoted the law of the random walk under $\omega$, or the annealed laws under $\P$ and $\hat{\P}$ respectively).
The new definitions below will still provide us with the random walk, but they will also contain the environment (unlike before).
To be more precise, let $(U_z)_{z\in\ZZ^{d+1}}$ be a collection of i.i.d.\ uniform random variables in $[0,1]$ defined on a probability space $(\Omega',\mathcal{T}',\P')$.
Let
\begin{equation}
\Omega_0 = \Omega \times \Omega', \quad \mathcal{T}_0 = \mathcal{T} \otimes \mathcal{T}', \text{ and } \bP_0 = \P \otimes \P'.
\end{equation}
We also set, for $\omega \in \Omega$,
\begin{equation}
  \P^\omega = \delta_{\{\omega\}} \otimes \P' \text{ and }
  \hat{\mathbb{P}} = \hat{\P} \otimes \P'.
\end{equation}

On these new space $(\Omega_0,\mathcal{T}_0)$, we can define the random walk $Z^z$ started at $z\in\ZZ^{d+1}$ driven by the environment $\omega$ inductively as follows:
\begin{equation}
  \label{e:def_Z}
  \left\{
    \begin{array}{l}
      Z^z_0=z; \\
      Z^z_{t+1}=Z^z_t+(g(\omega_{Z^z_t},U_{Z^z_t}),1), \quad t \in \NN.
    \end{array}\right.
\end{equation}
We simply write $Z$ for $Z^0$, where $0$ is the origin in $\ZZ^{d+1}$. Note that, with this construction, we do have equality \eqref{e:jump_law}.

Finally, define a shifting operator by setting, for any field $(F_y)_{y\in\ZZ^{d+1}}$ indexed by space-time points and $z,z_0\in\ZZ^{d+1}$,
\begin{equation}\label{e:translated}
 (F\circ\theta_{z_0})_{z}=F_{z_0+z}.
\end{equation}

\section{Independent decomposition}\label{s:ind_dec}

The field $\eta$ should be interpreted as a collection of amnesia traps for the random walk, i.e., places that allow it to renew and completely forget the past environment $\omega$ that it visited. However, it might be the case that the field $\eta$ itself still carries information of the environment $\omega$ into the future and keeps us from getting a decomposition of the walk into independent parts (notice that, while $\mathcal{F}_{0}^{-}$ contains the information about $\eta$ in the cone of the past $C_{0}^{-}$, the same is not true for $\mathcal{F}_{0}^{+}$ and the restriction $\eta|_{C_{0}^{+}}$, see Equation~\eqref{eq:sigma_algebras}). For this reason, we now introduce a resampling operator that deals with this undesired possibility. The existence of this resampling operator will be used below in Equation~\eqref{eq:measurable_iid}.

Recall that $\hat{\P}$ denotes the distribution $\P$ conditioned on $\eta_{0}=1$. Let us now define the properties that our resampling operator $\mathcal{R}$ should satisfy in order to solve our independence issues. It should use the information about $\omega$ in the cone of the future and some additional randomness to resample the value of the field $\eta$.

We will perform now the last enhancement of our probability space in order to include an additional collection of independent uniform random variables that will be used in combination with the resampling operator introduced above.
Let
\begin{equation}
\Omega_1 = \Omega \times \Omega' \times \Omega', \quad \mathcal{T}_1 = \mathcal{T} \otimes \mathcal{T}' \otimes \mathcal{T}', \text{ and } \bP = \P \otimes \P' \otimes \P'.
\end{equation}
We also set, for $\omega \in \Omega$, $\P^\omega = \delta_{\{\omega\}} \otimes \P'$ and $\hat{\mathbb{P}} = \hat{\P} \otimes \P'$ as before.
To distinguish the two collections of uniform random variables defined above we denote the first ones $U_z$, and the second $\tilde{U}_z$ for $z \in \mathbb{Z}^{d + 1}$.
The $U$'s will be used for the randomness of the random walk, while the $\tilde{U}$'s will be used for resampling with the help of the lemma below.

{\begin{definition}\label{d:resampling}
 $\mathcal R: S^{\Z^{d+1}}\times [0,1]\rightarrow \{0,1\}^{\Z^{d+1}}$ is a resampling operator under $\mathbb P$ if
 \begin{enumerate}
\item
$\big( \omega|_{C_0^{+}}, \mathcal{R}(\omega, \tilde{U}_0) \big) $ is a random variable with the same distribution as $ (\omega|_{C_0^{+}}, \eta)$ under $\mathbb P$,

\item There exists an $\big( \mathcal{F}_{0}^{+} \otimes \mathcal{B}[0,1] \big)$-measurable\footnote{Here, $\mathcal{B}[0,1]$ denotes the Borel $\sigma$-algebra on the interval $[0,1]$.} function $\tilde{\mathcal{R}}$ such that $\mathcal{R}(\omega, \tilde{U}_0) = \tilde{\mathcal{R}}(\omega, \tilde{U}_0)$ for $\mathbb P$-almost every $(\omega, \tilde{U})$.
\end{enumerate}
\end{definition}

Note that, if $\tilde{\mathcal R}$ is a resampling operator, since $\tilde{ \mathcal R}$ is $\big( \mathcal{F}_{0}^{+} \otimes \mathcal{B}[0,1] \big)$-measurable, $\tilde{\mathcal{R}}(\omega, \tilde{U}_0)$ above only depends on $\omega$ through its restriction to the cone of the future $C_0^+$.}
%

Recall the shift operator introduced in~\eqref{e:translated}. We write $\mathcal{R}_{z}(\omega, u) = \mathcal{R}(\omega \circ \theta_{z}, u) \circ \theta_{-z}$ for the resampling operator that uses the cone fixed at $z \in \Z^{d+1}$.

Let us note that, although it is not immediately clear that the operator above is always well-defined, this is a consequence of the existence of regular conditional probabilities, as stated in the next lemma.

\begin{lemma}\label{l:resampling_lemma}
Under the Random Markov Property, there exists a resampling operator under $\mathbb P$ as in Definition~\ref{d:resampling}.
\end{lemma}

We postpone the proof of this lemma to the end of this section and define now a sequence of random times for the walk starting at $z$.
More precisely, let the sequence $(T^{z}_{k})_{k \geqslant 0}$ as
\begin{equation}
  \label{e:T_k}
  \begin{array}{l}
    T^{z}_{0} = T^z_0(\omega) = \inf \big\{ t \geqslant 0,\, \eta_{Z^z_t} = 1 \big\}, \\
    T^z_{k + 1} = T^{z}_{k+1}(\omega) = \inf \Big\{ t > T^z_k,\, \mathcal{R}_{Z^z_{T^{z}_{k}}}\bigl(\omega, \tilde{U}_{Z^z_{T^{z}_{k}}}\bigr)_{Z^{z}_{t}} = 1 \Big\}, \text{ for all } k \in \NN,
  \end{array}
\end{equation}
where the infimum of the empty set is taken to be $+\infty$. We will often write simply $T_k$ without the upper script $z$ to refer to the times $T_k^0$ associated with the random walk starting at the origin. Observe that when $\eta_0 = 1$, we have $T_0 = 0$.

We also set, for $k\in\NN$,
\begin{equation}
  \label{e:Y_k}
  Y_k = \left\{\begin{array}{ll}
                 Z_{T_k} & \text{if $T_k<\infty$;}\\
                 \infty  & \text{otherwise.}
               \end{array}\right.
\end{equation}
We set $Y_{k+1}-Y_k= +\infty$ if $T_{k+1}= +\infty$.

The next lemma is quite general and does not require any decoupling assumption on $\eta$.

\begin{lemma}
  \label{l:indep_decomp}
    Let $n\geqslant 0$ and let $f_0,\ldots,f_n$ be bounded measurable functions on $\ZZ^{d+1}$. We have
    \begin{align}\label{e:big_equality}
     \mathbb{E}&\left[\mathbf{1}_{T_{n+1}<\infty}\,f_0(Y_1-Y_0)\cdots f_n(Y_{n+1}-Y_n)\,\, \middle| \,\eta_0=1\right]\nonumber\\
     &=\prod_{k=0}^n \bE[\mathbf{1}_{T_1<\infty}\,f_k(Y_1)\, | \,\eta_0=1].
    \end{align}
  In particular, under $\mathbb{P} \big( \cdot | \eta_0 = 1\big)$, if $T_1<\infty$ almost surely, then the $(Y_{k + 1} - Y_k)_{k\geqslant 0}$ are i.i.d. and have the same distribution as $Y_1=(X_{T_1}, T_1)$.
\end{lemma}

We will see in Lemma \ref{l:moments} that whenever $\alpha>0$, the $\alpha$-decoupling property implies that $T_1<\infty$ almost surely, so the second part of the statement will hold.
\begin{proof}
We first need to strengthen the Random Markov Property. Let
\begin{equation}
\mathcal{G}_z = \mathcal{F}_z^- \vee (U_y)_{y \in \mathcal{C}_z^- \setminus \{z\}}
\end{equation}
denote the $\sigma$-algebra generated by $\mathcal{F}_z^-$ and the collection of random variables $(U_y)_{y \in \mathcal{C}_z^- \setminus \{z\}}$.

Note that for any bounded measurable function $f$ on $(S\times [0,1]^2)^{\mathcal{C}_z^+}$, we have
\begin{align*}
\bE&\left[f\left((\omega_y,U_y,\tilde U_y)_{y\in\mathcal{C}_z^+}\right) \middle| \mathcal{G}_z\right] \eta_z\\
&=\bE\left[f\left((\omega_y,U_y,\tilde U_y)_{y\in\mathcal{C}_z^+}\right) \middle| \mathcal{F}_z^-\right] \eta_z\\
&=\int_{([0,1]^{2})^{\mathcal{C}_z^+}} \E\left[f\left((\omega_y,u_y,\tilde u_y)_{y\in\mathcal{C}_z^+}\right) \middle| \mathcal{F}_z^-\right] \eta_z \prod_{y\in\mathcal{C}_z^+} \mathrm{d}u_y \mathrm{d}\tilde u_y\\
&=\int_{([0,1]^{2})^{\mathcal{C}_z^+}} \E\left[f\left((\omega_y,u_y,\tilde u_y)_{y\in\mathcal{C}_z^+}\right) \middle| \eta_z\right] \eta_z \prod_{y\in\mathcal{C}_z^+} \mathrm{d}u_y \mathrm{d}\tilde u_y\\
&=\bE\left[f\left((\omega_y,U_y,\tilde U_y)_{y\in\mathcal{C}_z^+}\right) \middle| \eta_z\right]  \eta_z\\
&=\bE\left[f\left((\omega_y,U_y,\tilde U_y)_{y\in\mathcal{C}_z^+}\right) \middle| \eta_z=1\right]  \eta_z.
\end{align*}
The first equality uses that $(U_y)_{y\in\mathcal{C}_z^- \setminus \{z\}}$ is independent of $\mathcal{F}_z^-$ and $(\omega_y,U_y,\tilde U_y)_{y\in\mathcal{C}_z^+}$. In the second and last equalities, we used that $(U_y,\tilde U_y)_{y\in\mathcal{C}_z^+}$ is independent of $\omega$ and $\mathcal{F}_z^-$. In the third equality, the Random Markov property was used.

Let us now show an auxiliary identity. Let $f$ be a bounded measurable function on $\ZZ^{d+1}$ and $g$ a bounded measurable function on $(S\times [0,1]^2)^{\mathcal{C}_0^+}$. Under $\eta_0=1$, define $\tilde T=\inf\{t>0,\,\eta_{Z_t}(\omega)=1\}$ and $\tilde Y=Z_{\tilde T}$. {Notice that, although $(\tilde{Y}, \tilde{T})$ and $(Y_1, T_{1})$ may differ, they have the same law when conditioned on $\eta_0=1$, as the only difference is that the latter is defined using a resampled environment, which has the same law as $\eta$.} We have
\begin{equation}\label{eq:big_calculation}
\begin{split}
    \mathbb{E}&\left[\mathbf{1}_{T_1<\infty}\,f(Y_1)\,g\left((\omega,U,\tilde U)\circ \theta_{Y_1}\vert_{\mathcal{C}_0^+}\right) \middle| \eta_0=1\right]\\
    &=\mathbb{E}\left[\mathbf{1}_{\tilde T<\infty}\,f(\tilde Y)\,g\left((\omega,U,\tilde U)\circ \theta_{\tilde Y}\vert_{\mathcal{C}_0^+}\right)\middle| \eta_0=1\right]\\
    &=\sum_{z = (y, t) \in \ZZ^{d} \times \Z_{+}} f(z)\, \bE\left[\mathbf{1}_{\tilde T = t}\,\mathbf{1}_{\tilde Y=z} \,g\left((\omega_y,U_y,\tilde U_y)_{y\in\mathcal{C}_z^+}\right) \middle| \eta_0=1\right]\\
    &=\sum_{z = (y, t) \in \ZZ^{d} \times \Z_{+}} f(z)\, \bE\left[\mathbf{1}_{\tilde T = t}\,\mathbf{1}_{\tilde Y=z} \,\bE\left[g\left((\omega_y,U_y,\tilde U_y)_{y\in\mathcal{C}_z^+}\right) \middle| \mathcal{G}_z\right] \middle| \eta_0=1\right]\\
    &=\sum_{z = (y, t) \in \ZZ^{d} \times \Z_{+}} f(z)\, \bE\left[\mathbf{1}_{\tilde T = t}\,\mathbf{1}_{\tilde Y=z} \,\bE\left[g\left((\omega_y,U_y,\tilde U_y)_{y\in\mathcal{C}_z^+}\right) \middle| \eta_z=1\right] \middle| \eta_0=1\right]\\
    &=\bE\left[\mathbf{1}_{\tilde T<\infty}\,f(\tilde Y) \middle| \eta_0=1\right]\,\bE\left[g\left((\omega_y,U_y,\tilde U_y)_{y\in\mathcal{C}_0^+}\right) \middle| \eta_0=1\right]\\
    &=\bE\left[\mathbf{1}_{T_1<\infty} \,f(Y_1) \middle| \eta_0=1\right]\,\bE\left[g\left((\omega,U,\tilde U)\vert_{\mathcal{C}_0^+}\right) \middle| \eta_0=1\right].
\end{split}
\end{equation}
In the first and last equalities, the equality in distribution between $(\tilde{Y}, \tilde{T})$ and $(Y_1, T_{1})$ is used. The third equality relies on the fact that, { for $z = (y, t) \in \ZZ^{d} \times \Z_{+}$, the events $\{\eta_0=1\}$, $\{ \tilde T = t\}$, and $\{ \tilde Y = z\}$} are $\mathcal{G}_z$-measurable. In the fourth equality, we used the strengthened Random Markov property, while the fifth equality follows from translation invariance.

We now verify~\eqref{e:big_equality} by induction in $n$. It is obvious when $n=0$. Suppose that~\eqref{e:big_equality} is satisfied for $n-1$ with some $n \geqslant 1$. Let $f_0, \ldots, f_n$ be bounded measurable functions on $\ZZ^{d+1}$ and notice that there exists a bounded measurable function $g$ as in~\eqref{eq:big_calculation} such that
\begin{equation}\label{eq:measurable_iid}
\mathbf{1}_{T_{n+1}<\infty}\,f_1(Y_2-Y_1)\cdots f_n(Y_{n+1}-Y_n)=g\left((\omega,U,\tilde U)\circ \theta_{Y_1}\vert_{\mathcal{C}_0^+}\right).
\end{equation}
This is where the resampling operator comes in handy. Indeed, here we are using that the sequence $(Y_{k+1}-Y_{k})_{k \geq 1}$ is measurable with respect to the environment in the cone of the future at $Y_{1}$, which is due to the fact that the resampling uses only information about the future environment in order to redefine the field $\eta$.

Applying~\eqref{eq:big_calculation} with $f=f_0$ and $g$, we get the result after noticing that under $\eta_0=1$, we have
$$g\left((\omega,U,\tilde U)\vert_{\mathcal{C}_0^+}\right)=\mathbf{1}_{T_n<\infty}\, f_1(Y_1)\cdots f_n(Y_n-Y_{n-1})$$ and applying the induction assumption.

Assume now that $\bP(T_1<\infty\,\vert\,\eta_0=1)=1$ and apply~\eqref{e:big_equality} with $n \in \NN$ and $f_0 = \cdots = f_n = 1$ to get $\bP(T_{n+1}<\infty\vert\,\eta_0=1)=1$. Therefore, \eqref{e:big_equality} becomes 
    \begin{align}\label{e:big_equality2}
     \mathbb{E}\left[f_0(Y_1-Y_0)\cdots f_n(Y_{n+1}-Y_n)\,\vert\,\eta_0=1\right]=\prod_{k=0}^n \bE[f_k(Y_1)\,\vert\,\eta_0=1],
    \end{align}
for any bounded measurable functions $f_0,\ldots,f_n$. This concludes the proof of lemma.
\end{proof}

We end this section with the proof of Lemma~\ref{l:resampling_lemma}.
\begin{proof}[Proof of Lemma~\ref{l:resampling_lemma}]
Recall that $\mathcal{F}^{+}_{0} = \sigma \big (\omega_{z}: z \in C^{+}_{0} \big)$. Since $\omega$ takes values in a Polish space, the regular conditional probability $K(\omega, \cdot) = \hat{\P} ( \, \cdot \, | \mathcal{F}^{+}_{0})(\omega)$ is well defined\footnote{This means that $K(\omega, A)$ is a transition kernel that satisfies: 1) $A \mapsto K(\omega, A)$ is a probability measure on $\mathcal{T}$ for $\hat{\P}$-almost every $\omega$. 2) $\omega \mapsto K(\omega, A)$ is measurable for every measurable $A$. 3) For every $A$, $\omega \mapsto K(\omega, A)$ is a $\hat{\P}$-version of $\hat{\E}(1_A| \mathcal{F}^{+}_{0})$.} for $\hat{\P}$-almost every $\omega$, {(see for example~\cite[Theorem 4.1.18]{durrett_book})}.

Although {the} kernel $K$ allows us to sample pairs $(\tilde{\omega}, \tilde{\eta})$, we only need the latter and use the variable {$\tilde{U}_0$  to sample\footnote{The fact that we can do that is a standard property of real-valued variables, but here we want to sample a whole field of real random variables. To do this, we may first use $\tilde{U}_0$ to generate an i.i.d.\ field of uniform random variables $(\tilde U_{0,z})_{z\in\mathbb Z^{d+1}}$. This can be done by writing the decomposition of $\tilde U_0$ into negative powers of $2$ (yielding i.i.d.\ $\rm{Ber}(1/2)$ variables indexed by $\mathbb N$), then using a bijection between $\mathbb N$ and $\Z^{d+1}\times\mathbb N$ to reindex them. The $\mathbb N$--indexed collection corresponding to $z\in\mathbb Z^{d+1}$ gives the decomposition of $\tilde U_{O,z}$ into negative powers of $2$. Then we fix an enumeration $z_1,z_2,\ldots$ of $\Z^{d+1}$ and use $U_{0,z_i}$ to sample $\tilde\eta_{z_i}$, conditionally on $\eta_{z_1},\ldots,\eta_{z_{i-1}}$.} a field $\tilde{\eta} \in \{0,1\}^{\Z^{d+1}}$ from $K(\omega, \cdot)$, thus defining $\mathcal{R}(\omega, \tilde{U}_0) = \tilde{\eta}$. }


For the first property of the resampling operator, simply notice that
{\begin{equation}
  \begin{split}
    \hat{\P} \big( \tilde{\eta} \in A, \omega|_{C_0^{+}} \in B \big) & = \hat{\E} \big[ \textbf{1}_{\omega|_{C_0^{+}} \in B} K(\omega, \{\eta \in A \}) \big] \\
    & = \hat{\E} \big[ \textbf{1}_{\omega|_{C_0^{+}} \in B}  \hat{\P} \big( \eta \in A \big| \mathcal{F}^{+}_{0} \big) \big] = \hat{\P} \big( \eta \in A, \omega|_{C_0^{+}} \in B \big),
  \end{split}
\end{equation}}
for any measurable sets $A \subset \{0,1\}^{\Z^{d+1}}$ { and $B \subset S^{C_0^{+}}$}.

Let us now verify the second property. Notice that, from the definition of the kernel $K(\cdot, \cdot)$, for every finite-dimensional set $A \subset \{0,1\}^{\Z^{d+1}}$, there exists a set $G_{A} \subset \Omega$ with $\hat{\P}(G_{A})=1$ such that $K(\omega, \{\eta \in A\}) = \hat{\P} (\{ \eta \in A\} | \mathcal{F}^{+}_{0})$ in $G_{A}$. Consider now $G = \cap_{A} G_{A}$ (the intersection runs over all finite-dimensional sets) and notice that an application of Dynkin's $\pi-\lambda$ Theorem yields that, for any measurable set $B \subset \{0,1\}^{\Z^{d+1}}$, $K(\omega, \{\eta \in B\}) = \hat{\P} ( \eta \in B | \mathcal{F}^{+}_{0})$ in $G$. By further restricting the set $G$, we may assume that, for all $\omega \in G$, $B \mapsto K(\omega, B)$ is a probability measure. In particular, in $G$, the function $\omega \mapsto K(\omega, \cdot)$, that maps configurations to probability measures, is measurable with respect to $\mathcal{F}^{+}_{0}$. The last ingredient missing is to note that, given $\tilde{U}$, the sampling of $\mathcal{R}$ is measurable with respect to the regular conditional probability kernel $K(\omega, \cdot)$. In particular, in the set  $G$, $\mathcal{R}(\omega, \tilde{U})$ coincides with an $\big( \mathcal{F}^{+}_{0} \otimes \mathcal{B}[0,1] \big)$-measurable function $\hat{\P}$-almost surely, concluding the verification of the second property.
\end{proof}

\section{Renewal tails}\label{s:renewal}

\subsection{Threatened points}

Throughout this section we consider the random walk $Z^z$ on top of the random environment $\omega$, governed by the law $\P^\omega$ defined in~\eqref{e:jump_law}.

We fix a collection of \emph{traps} $\Sigma \subseteq \mathbb{Z}^{d + 1}$ and estimate the time it takes for the random walk to hit one of these traps.
More precisely, consider the stopping time
\begin{equation}
  \label{e:T_trap}
  T^z = \inf\{ t>0; Z^z_t \in \Sigma \}.
\end{equation}
The dependence of quantities such as $T^z$ on the set $\Sigma$ will be omitted through the section, since $\Sigma$ is considered fixed.

Later, when we are going to apply the results of this section to our random walk, the set $\Sigma$ will be chosen as $\{ z \in \mathbb{Z}^{d + 1}; \eta_z = 1 \}$, so that conditioned on $\eta_0=1$, the stopping time $T^z$ will coincide with $T^z_1$ introduced in \eqref{e:T_k}.
However, the results of this section work for any set $\Sigma \subseteq \mathbb{Z}^{d + 1}$ and they are not constrained to being given by the times where we can construct a decoupling of $\omega$'s past and future (as in the definition of $\eta$).

Our goal is to give conditions under which the random walk falls on some trap sufficiently quickly. In the definition below, we denote by $\pi_{d+1}: \Z^{d+1} \to \Z$ the projection in the $(d+1)$-coordinate and by $\bar{\pi}_{d}: \Z^{d+1} \to \Z^{d}$ the projection in the first $d$ coordinates, so that a point in $\Z^{d+1}$ can be written as $x = (\bar{\pi}_{d}(x), \pi_{d+1}(x)) \in \Z^{d} \times \Z$.
\begin{definition}
  Let $a < b \in \NN$.
  We say that a path $\sigma: \llbracket a, b \rrbracket \rightarrow \ZZ^{d + 1}$ is an $R$-allowed path if
  \begin{itemize}
  \item $\pi_{d+1}(\sigma(s)) = s$, for every $s\in\llbracket a,b\rrbracket$: ``$d + 1$ is time'';
  \item $\bar{\pi}_{d} \circ \sigma$ is $R$-Lipschitz in the $|\cdot|_\infty$ norm, meaning that
    \begin{equation}
      |\bar{\pi}_{d} \circ \sigma(x) - \bar{\pi}_{d} \circ \sigma(y)|_{\infty} \leq R |x-y|, \text{ for all } x,y \in \llbracket a, b \rrbracket.
    \end{equation}
  \end{itemize}
  We call $b-a$ the length of the allowed path $\sigma$. By construction, the sample paths of the random walks $Z^z$ on any time interval are $R$-allowed paths.
\end{definition}

\begin{definition}
  \label{d:threat}
  For $H \geqslant 1$, we say that a point $z = (x, t)\in\ZZ^{d+1}$ is $H$-threatened (with respect to $\Sigma$) if there exists a $1$-allowed path $\sigma: \llbracket t, t + H \rrbracket\to \ZZ^{d+1}$ such that $\sigma(t) = z$ and $\sigma$ passes through some point $z' \in \Sigma$.
\end{definition}

Observe that we require the points $z$ and $z'$ to be connected through a $1$-allowed path instead of an $R$-allowed path.
The reason behind this is that we intend to use the uniform ellipticity provided by \eqref{e:ellipticity}, which guarantees that our random walk has a positive probability to follow the exact steps of any $1$-allowed path.
Therefore, whenever $z=(x,t)$ is $H$-threatened with respect to $\Sigma$, we have, for every $\omega$,
\begin{equation}\label{e:consequence_threat}
  \P^\omega \left( \text{there exists $s\in\llbracket t, t+H\rrbracket$ such that $Z_{s}^z \in \Sigma$} \right) \geqslant \uc{c:ellipticity}^H,
\end{equation}
which can be deduced from \eqref{e:ellipticity} by instructing the random walk to follow step-by-step the $1$-allowed path from Definition~\ref{d:threat}.

We now define the \emph{number of threats} along a given path.
\begin{definition}
  \label{d:num_threats}
  Given an $R$-allowed path $\sigma: \llbracket a, b \rrbracket\to \ZZ^{d+1}$ and $H \geqslant 1$, let
  \begin{equation}
    \label{e:num_threats}
    M(\sigma, H) := \Big| \big\{ t \in \llbracket a, b \rrbracket : \sigma(t) \text{ is $H$-threatened and $\pi_{d+1}(\sigma(t)) \in H \mathbb{Z}$} \big\} \Big|.
  \end{equation}
\end{definition}

Observe that we only count a threat if it is encountered at a time which is a multiple of $H$, this is done in order to keep the attempts to fall into different traps independent.

The connection between threats and the desired estimates on the tail of $T^z$ are made evident in the following result.
\begin{lemma}
  \label{l:fall_on_trap}
  For every $\omega\in\Omega$ and $H,J \in \NN^*$,
  \begin{equation}
    \label{e:fall_on_trap}
     \P^\omega \left( T^0 > J H \right) \leqslant (1 - \uc{c:ellipticity}^H)^{M_J},
  \end{equation}
  where $M_J = M_J(\Sigma)$ is given by
  \begin{equation*}
    M_J := \inf \Big\{ M(\sigma, H); \sigma: \llbracket 0, J H \rrbracket \rightarrow \ZZ^{d+1} \text{ $R$-allowed path starting at $0$} \Big\}.
  \end{equation*}
\end{lemma}

Before we provide the proof of the lemma, let us explain the simple intuition behind it.
Note that the path $\sigma'$ given by the random walk $(Z^0_t)_{t \in \llbracket 0, J H \rrbracket}$ itself is an \mbox{$R$-allowed} path starting at zero, therefore $M(\sigma', H) \geqslant M_J$.
In other words, independently of what the random walk does, it will be threatened $M_J$ times and in well separated locations.
Also, in each of these threatened times, by \eqref{e:consequence_threat}, there is a chance of at least $\uc{c:ellipticity}^H$ that the walk falls into the corresponding trap.
The proof below makes the above intuition rigorous.

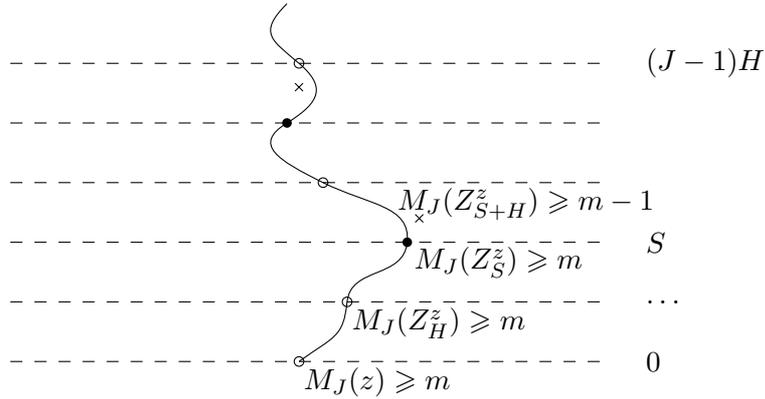
\begin{figure}
  \centering
  \begin{tikzpicture}[x=0.75pt,y=0.75pt,yscale=-60,xscale=60]
    \draw  [dash pattern={on 4.5pt off 4.5pt}]  (1,3.0) -- (6,3.0) ;
    \draw  [dash pattern={on 4.5pt off 4.5pt}]  (1,0.5) -- (6,0.5) ;
    \draw  [dash pattern={on 4.5pt off 4.5pt}]  (1,1.0) -- (6,1.0) ;
    \draw  [dash pattern={on 4.5pt off 4.5pt}]  (1,1.5) -- (6,1.5) ;
    \draw  [dash pattern={on 4.5pt off 4.5pt}]  (1,2.0) -- (6,2.0) ;
    \draw  [dash pattern={on 4.5pt off 4.5pt}]  (1,2.5) -- (6,2.5) ;
    \draw (6.2,3.0) node [anchor=west] {$0$};
    \draw (6.2,2.5) node [anchor=west] {$\dots$};
    \draw (6.2,2.0) node [anchor=west] {$S$};
    \draw (6.2,0.5) node [anchor=west] {$(J-1) H$};
    \draw  (3.40,3.00) .. controls (3.69,2.79) and (3.77,2.80) .. (3.80,2.50) .. controls (3.83,2.20) and (4.27,2.31) .. (4.30,2.00) .. controls (4.33,1.69) and (3.98,1.69) .. (3.60,1.50) .. controls (3.22,1.31) and (3.01,1.18) .. (3.30,1.00) .. controls (3.59,0.82) and (3.62,0.69) .. (3.40,0.50) .. controls (3.18,0.31) and (3.05,0.25) .. (3.30,0) ;
    \draw (3.40, 3) circle (0.04);
    \draw (3.80, 2.5) circle (0.04);
    \fill (4.3, 2) circle (0.04);
    \draw (4.4, 1.8) node[cross] {};
    \draw (3.60, 1.5) circle (0.04);
    \fill (3.30, 1) circle (0.04);
    \draw (3.4, 0.7) node[cross] {};
    \draw (3.4, .5) circle (0.04);
    \draw (3.42,3.03) node [anchor=north west][inner sep=0.75pt] {$M_{J}(z) \geqslant m$};
    \draw (3.82,2.53) node [anchor=north west][inner sep=0.75pt] {$M_{J}(Z_H^{z}) \geqslant m$};
    \draw (4.34,2.02) node [anchor=north west][inner sep=0.75pt]  {$M_{J}(Z_S^z) \geqslant m$};
    \draw (4.2,1.52) node [anchor=north west][inner sep=0.75pt]  {$M_{J}(Z_{S + H}^z) \geqslant m - 1$};
  \end{tikzpicture}
  \caption{An illustration of the proof of Lemma~\ref{l:fall_on_trap}.
    The horizontal dashed lines represent times that are multiple of $H$.
    The filled dots along the random walk trajectory represent threatened points, while the crosses represent traps in $\Sigma$.
    Note how the lower bound on $M_J$ drops from $m$ to $m - 1$ after crossing the first threat.}
  \label{f:fall_on_trap}
\end{figure}

\begin{proof}
  We follow an induction argument on $M_J$, but in order to set this up, we need to introduce a more general version of $M_J$ that includes alternative starting points for the path.
  Given $z = (x, j H) \in \mathbb{Z}^d \times (H \mathbb{Z})$ with $j \in \llbracket 0, J - 1 \rrbracket$, define
  \begin{equation}
    \label{e:M_J_z}
    M_J(z) = \inf \Big\{ M(\sigma, H); \sigma: \llbracket jH, JH \rrbracket \to \Z^{d+1} \text{ starts at $z$ and is $R$-allowed} \Big\},
  \end{equation}
  which depends solely on $\Sigma$.
  Recall the definition of $T^{z}$ from~\eqref{e:T_trap}. Our objective is to prove that, for all $\Sigma$,
  \begin{equation}
    \label{e:fall_on_trap_general}
    \begin{gathered}
      \text{for every $z = (x, j H)$ (with $j = 0, \dots, J - 1$), for every $\omega$ and $m \geqslant 1$}\\
      \text{if $M_J(z) \geqslant m$, then } \P^\omega(T^z > (J - j) H) \leqslant (1 - \uc{c:ellipticity}^H)^{m}.
    \end{gathered}
  \end{equation}
  Note that \eqref{e:fall_on_trap_general} is enough to establish Lemma~\ref{l:fall_on_trap} if we apply it to $z = (0, 0)$ and $m = M_J$.
  We now turn to the proof of \eqref{e:fall_on_trap_general} by using induction on $m$.

  It is clear that the claim is valid for $m = 0$ and from now on we suppose that $m > 0$ and \eqref{e:fall_on_trap_general} is valid for $m - 1$.
  Take then an arbitrary point $z = (x, j H)$ with $j \in \{0, \dots, J-1\}$ such that $M_J(z) \geqslant m$ and we turn to the proof of \eqref{e:fall_on_trap_general}.
  Define first the stopping time (with respect to the filtration $\sigma( Z^z_s, s \leqslant t)$ on the space of random walk trajectories)
  \begin{equation}
    S = \inf \big\{ s \in \{0, H, 2H, \dots\}; Z^z_{s} \text{ is $H$-threatened} \big\},
  \end{equation}
  which is smaller than or equal to $(J - j - 1) H$ almost surely because $M_J(z) \geqslant m \geqslant 1$.
  We then estimate, using the Strong Markov Property twice (at times $S$ and $S+H$, see Figure~\ref{f:fall_on_trap}),
  \begin{equation}
    \label{e:inductive_bound_MJ}
    \begin{split}
      \P^\omega (T^z > & (J - j) H)
      =\E^\omega \Big[ \P^\omega \left(T^z > (J - j) H \big| Z^z_0, \dots, Z^z_S\right) \Big]\\
      & \leqslant \E^\omega \Big[ \P^\omega \left( T^{Z^z_S} > (J - j) H - S \right) \Big] \\
      & \leqslant \E^\omega \Big[ \P^\omega \left( T^{Z^z_S} > H, T^{Z^z_{S + H}} > (J - j - 1) H - S \right) \Big] \\
    & \leqslant \E^\omega \bigg[ \Big( 1 - \uc{c:ellipticity}^H \Big) \sup_{z' \text{ ; $M_J(z') \geqslant m-1$}} \P^\omega \left( T^{z'} > (J - j-1)H - S \right) \bigg] \\
      & \leqslant  \Big( 1 - \uc{c:ellipticity}^H \Big)^m.
    \end{split}
  \end{equation}
  Notice that, in the estimate above, the supremum is taken over all points $z' = (x, (j+1)H+S)$ such that $M_J(z') \geqslant m-1$, which allows for the direct application of the induction hypothesis. This concludes the proof of \eqref{e:fall_on_trap_general} and consequently the lemma.
\end{proof}

\subsection{Finding traps on intervals}\label{ss:traps}

This section is dedicated to proving that long intervals are very likely to intersect a trap.
This is done using a renormalization technique that is inspired by \cite{cavalinhos}.
We start defining the sequence of scales
\begin{equation}
  \label{e:scales}
  L_k := 4^k, \qquad \ell_k = \Big\lfloor \frac{L_k}{2 d k^2} \Big\rfloor.
\end{equation}

Given integers $a < b$, we introduce the event
\begin{equation}
  \label{e:Fab}
  F(a, b) := \big\{ \eta_{(0, a)}
  = \eta_{(0, a + 1)} = \cdots = \eta_{(0, b - 1)} = 0 \big\}.
\end{equation}

{Recall, as already mentioned immediately below~\eqref{e:T_trap}, that the set of traps we consider here is precisely $\Sigma = \{ z \in \Z^{d+1}: \eta_{z}=1\}$.} We want to show that
\begin{equation}
  \label{e:D_k}
  q_k := \P\big( F ( 0, L_k ) \big),
\end{equation}
decays fast enough. To that end, we first establish a recursive inequality between $q_{k+1}$ and $q_k$, \eqref{e:D_induction}. This inequality has the desirable property that, if $q_k\leq cL_k^{-\alpha}$ for some large $k$, the estimate propagates to higher $k$'s; this induction step is the proof of Lemma~\ref{l:rod}. The missing ingredient is then a triggering property, which consists in finding a suitable $k$ to initiate the recursion; this is done with the help of Lemma~\ref{l:rod_weak}. 

Observe that
\begin{equation}
  \label{e:F_cascades}
  F \big( 0, L_{k + 1} \big)
  \subseteq F \big( 0, L_k \big) \cap F \big( 3 \cdot L_k, L_{k + 1} \big).
\end{equation}

From \eqref{e:F_cascades}, applying the $\alpha$--decoupling property (recall Definition \ref{d:decoupling}) and translation invariance,

\begin{equation}
  \label{e:D_induction}
  q_{k + 1} \leqslant q_k^2 + \uc{c:decoupling} L_k^{-\alpha}.
\end{equation}

Recall the constant $\uc{c:eta_zero}$ from \eqref{e:eta_zero}.
The next lemma shows that $q_k$ decays with $k$.
\nc{c:decoupling_decay}
\nc{c:decoupling_decay_2}
\begin{lemma}
  \label{l:rod_weak}
  Under the $\alpha$--decoupling property, there exist $\uc{c:decoupling_decay} = \uc{c:decoupling_decay}(\alpha, \uc{c:eta_zero}, \uc{c:decoupling}) \in (0, 1)$ and $\uc{c:decoupling_decay_2} = \uc{c:decoupling_decay_2}(\alpha, \uc{c:eta_zero}, \uc{c:decoupling}) {> 0}$ such that
  \begin{equation}
    \label{e:rod_weak}
    q_k \leqslant \uc{c:decoupling_decay_2} \uc{c:decoupling_decay}^k, \text{ for every $k \geqslant 0$}.
  \end{equation}
\end{lemma}

\begin{proof}
  We first fix a constant $\uc{c:decoupling_decay} = \uc{c:decoupling_decay}(\alpha, \uc{c:eta_zero}) \in (0, 1)$ such that
  \begin{equation}
    \label{e:boot_weak}
    \uc{c:decoupling_decay} > 4^{-\alpha/2} \qquad \text{and} \qquad
    \uc{c:decoupling_decay}^2 \geqslant \P\big( \eta_{(0, 0)} = 0 \big).
  \end{equation}
Note that this is possible using \eqref{e:eta_zero}. Choose also an integer $k_0 = k_0(\alpha, \uc{c:eta_zero}, \uc{c:decoupling}) \geqslant 1$ such that
  \begin{equation}
    \label{e:k_zero_rod_weak}
    \uc{c:decoupling} 4^{- \alpha \big( k_0 - \frac{1}{2} \big)} \leqslant 1 - \uc{c:decoupling_decay}.
  \end{equation}

  We want to prove by induction in $k \geqslant 2$ that
  \begin{equation}
    \label{e:q_k_decay}
    q_{k_0 + k} \leqslant \uc{c:decoupling_decay}^k,
  \end{equation}
  which is sufficient to establish \eqref{e:rod_weak}.

  Taking $k = 2$, observe that
  \begin{equation}
    q_{k_0 + 2} =\P\big( F(0, L_{k_0 + 2}) \big)
    \leqslant \P\big( \eta(0, 0) = 0 \big)
    \overset{\eqref{e:boot_weak}}\leqslant \uc{c:decoupling_decay}^2,
  \end{equation}
  so that \eqref{e:q_k_decay} holds for $k = 2$.

  Assuming now the validity of~\eqref{e:q_k_decay} for a certain value of $k \geqslant 2$, we get
  \begin{equation}
    \begin{split}
      \frac{q_{k_0 + k + 1}}{\uc{c:decoupling_decay}^{k+1}}
      & \overset{\eqref{e:D_induction}}\leq
      \frac{q_{k_0 + k}^2 + \uc{c:decoupling} L_{k_0 + k}^{-\alpha}}{\uc{c:decoupling_decay}^{k+1}}
      \overset{\eqref{e:q_k_decay}}\leq
      \uc{c:decoupling_decay}^{k - 1} + \uc{c:decoupling} \frac{L_{k_0 + k}^{-\alpha}}
      {\uc{c:decoupling_decay}^{k + 1}}\\
      & \overset{\eqref{e:boot_weak}}\leqslant \uc{c:decoupling_decay}^{k - 1} + \uc{c:decoupling} 4^{-\alpha(k_0 + k) + \alpha(k+1)/2} \overset{\eqref{e:k_zero_rod_weak}} \leqslant 1.
    \end{split}
  \end{equation}
  This finishes the proof of \eqref{e:q_k_decay} and we can bootstrap it into \eqref{e:rod_weak} by choosing $\uc{c:decoupling_decay}$ sufficiently close to one.
\end{proof}

Note that the decay obtained from the above lemma is slow, this will be improved in the next lemma.

\nc{c:rod}
\begin{lemma}
  \label{l:rod}
  Under the $\alpha$--decoupling property, there exists $\uc{c:rod} = \uc{c:rod}(\alpha, \uc{c:decoupling}, \uc{c:eta_zero})$ such that
  \begin{equation}
    \label{e:rod}
    q_k \leqslant \uc{c:rod} L_k^{-\alpha},
  \end{equation}
  for every $k \geqslant 0$.
\end{lemma}

\begin{proof}
  Using Lemma~\ref{l:rod_weak}, we can obtain $k_0 = k_0(\alpha, \uc{c:decoupling}, \uc{c:eta_zero}) \geqslant 1$ such that
  \begin{equation}
    \label{e:k_zero_rod}
    q_{k_0 + 1} \leqslant \frac{1}{2} L_1^{-\alpha} \qquad \text{and} \qquad 2 \uc{c:decoupling} 4^{-\alpha(k_0 - 1)} \leqslant \frac{1}{2}.
  \end{equation}

  We are going to prove by induction that,
  \begin{equation}
    \label{e:q_k_decay_2}
    q_{k_0 + k} \leqslant \frac{1}{2} L_k^{-\alpha}, \text{ for every $k \geqslant 1$}.
  \end{equation}
  Note that the case $k = 1$ is covered in the left hand side of \eqref{e:k_zero_rod}. Suppose that \eqref{e:q_k_decay_2} holds for a certain $k \geqslant 1$ and estimate
  \begin{equation}
    \begin{split}
      \frac{q_{k_0 + k + 1}}{(1/2)L_{k + 1}^{-\alpha}}
      & \overset{\eqref{e:D_induction}}\leqslant 2 \frac{q_{k_0 + k}^2
        + \uc{c:decoupling} L_{k_0 + k}^{-\alpha}}{L_{k + 1}^{-\alpha}}
        \overset{\eqref{e:q_k_decay_2}}\leqslant 2 \cdot \frac{1}{4} \cdot L_k^{-2 \alpha} L_{k + 1}^\alpha
        + 2 \uc{c:decoupling} L_{k + k_0}^{-\alpha} L_{k + 1}^{\alpha}\\
      & \leqslant \frac{1}{2} 4^{-\alpha (k - 1)}
        + 2 \uc{c:decoupling} 4^{-\alpha (k_0 - 1)}
        \overset{\eqref{e:k_zero_rod}}\leqslant \frac{1}{2} + \frac{1}{2} = 1.
    \end{split}
  \end{equation}
  This finishes the proof of \eqref{e:q_k_decay_2} by induction.
  Estimate~\eqref{e:rod} now follows from \eqref{e:q_k_decay_2} by properly choosing $\uc{c:rod} = \uc{c:rod}(k_0) = \uc{c:rod}(\alpha, \uc{c:decoupling}, \uc{c:eta_zero})$ to adjust the exponent for $k > k_{0}$ and to cover the remaining cases.
\end{proof}

\begin{corollary}\label{cor:rodL}
	\nc{c:rodL}
	There exists $\uc{c:rodL}=\uc{c:rodL}(\alpha, \uc{c:decoupling}, \uc{c:eta_zero})$ such that, for any positive integer $L$
	\begin{equation}
		\P(F(0,L))\leqslant \uc{c:rodL}L^{-\alpha}.
		\end{equation}
	\end{corollary}
\begin{proof}
	The proof is a simple interpolation: fix $L$ and $k$ such that $L_{k} \leqslant L < L_{k+1}$ and note that
	\begin{align}
		\P(F(0,L)) \leqslant q_k \leqslant \uc{c:rod} L_k^{-\alpha} \leqslant \uc{c:rod}4^\alpha L_{k+1}^{-\alpha} \leqslant \uc{c:rod}4^\alpha L^{-\alpha},
        \end{align}
        proving the result.
	\end{proof}

With the uniform ellipticity~\eqref{e:ellipticity}, this implies that the probability that \emph{every} point in a large box of $\Z^d$ be threatened is large (recall Definition~\ref{d:threat}). This will be useful in the next section where we show that every allowed path meets many threatened points with high probability.

\begin{corollary}\label{cor:threatsinbox}
	\nc{c:threatsinbox}
	There exists $\uc{c:threatsinbox}=\uc{c:threatsinbox}(\alpha, \uc{c:decoupling}, \uc{c:eta_zero})$ such that for any $w\in\Z^{d+1}$, for any integer $L$,
	\begin{equation}
		\P(\exists z\in w+\llbracket 0,L-1\rrbracket^d\times\{0\} \text{ s.t. $z$ is not $L$--threatened})\leqslant \uc{c:threatsinbox}L^{-\alpha}.
		\end{equation}
	\end{corollary}
\begin{proof}
	If $L$ is odd, let us call $x_0=(\lfloor L/2\rfloor,\ldots,\lfloor L/2\rfloor )\in\Z^d$ the middle point of $\llbracket 0,L-1\rrbracket^d$. If there is a trap on $w+\{x_0\}\times\llbracket \lfloor L/2\rfloor,L\rrbracket$, by construction, every point in $w+\llbracket 0,L-1\rrbracket^d\times\{0\}$ is $L$--threatened. We then use translation invariance and apply Corollary~\ref{cor:rodL}. If $L$ is even, we have, using translation invariance and a union bound,
	\begin{align*}
	 \P&(\exists z\in w+\llbracket 0,L-1\rrbracket^d\times\{0\} \text{ s.t. $z$ is not $L$--threatened})\\
	 &\leqslant \P(\exists z\in w+\llbracket 0,L-1\rrbracket^d\times\{0\} \text{ s.t. $z$ is not $(L-1)$--threatened})\\
	 &\leqslant 2^d \P(\exists z\in w+\llbracket 0,L-2\rrbracket^d\times\{0\} \text{ s.t. $z$ is not $(L-1)$--threatened})\\
	 &\leqslant 2^d{ \uc{c:rodL}}(L-1)^{-\alpha}.
	\end{align*}
    The result follows by choosing $\uc{c:threatsinbox}$ accordingly.
\end{proof}

\subsection{Number of threats on allowed paths}

Having established that, with high probability, one can find threatened points on segments, this section moves on to prove that it is very likely that every allowed path finds many traps along its way.

For $k\in\NN$, $H\geqslant 1$ and $w\in\ZZ^{d+1}$, we introduce the box and its centered base:
\begin{align*}
 &B_{k}=\llbracket -RL_k,(R+1)L_k-1 \rrbracket^d\times\llbracket 0,L_k-1\rrbracket\subseteq \ZZ^{d+1};\\
 &R_{k}=\llbracket 0,L_k-1\rrbracket^d\times\{0\}.
\end{align*}
We also define the event that an allowed path with too few threats:
\begin{equation}
A_{k,H}(w)=\bigg\{ \begin{array}{c}
\text{there exists an $R$--allowed path } \sigma \text{ of length } L_k-1 \\
 \text{ starting in } w+R_k \text{ such that } M(\sigma,H)< k^2
 \end{array} \bigg\},
\end{equation}
and $A_{k,H}=A_{k,H}(0)$.

The goal of this section is to show that the above event is unlikely:
\begin{lemma}\label{l:estimate_A}
For every $\beta<\alpha$, there exist $\uc{c:triggering}=\uc{c:triggering}(\alpha,\beta,\uc{c:decoupling}, \uc{c:eta_zero})>0$ and two integers $k_0=k_0(\alpha,\beta,\uc{c:decoupling}, \uc{c:eta_zero})$ and $H=H(\alpha,\beta,\uc{c:decoupling}, \uc{c:eta_zero})$ satisfying, for every $k\geqslant k_0$,
\begin{equation}\label{e:estimate_A}\P(A_{k,H})\leqslant \uc{c:triggering} L_k^{-\beta}.\end{equation}
\end{lemma}

The proof of this lemma relies on a renormalization strategy, similar to what we did in Section~\ref{ss:traps}, but slightly more involved (recall the discussion after \eqref{e:D_k}).

\subsubsection{Renormalization argument}
\label{sss:renormalization}
For now, $H$ is just a fixed parameter; it will be chosen in Section \ref{sss:triggering}. In this section, we establish a recursive inequality between $\bP(A_{k+1,H})$ and $\bP(A_{k,H})$.


The core of the argument is the following. Fix $k\geqslant 3$, $H\in\NN^*$ and suppose that $A_{k+1,H}$ occurs. We are therefore given an $R$--allowed path $\sigma:\llbracket 0,L_{k+1}-1\rrbracket\rightarrow \ZZ^{d+1}$ starting in $R_{k+1}$ such that $M(\sigma,H)<(k+1)^2$. Recalling \eqref{e:num_threats}, this means that there are at most $(k+1)^2$ times $t\in \left\{0,H,\ldots,\left\lfloor \frac{L_{k+1}-1}{H}\right\rfloor H\right\}$ such that $\sigma(t)$ is $H$-threatened.

Since $\sigma$ is an $R$--allowed path starting in $R_{k+1}$, by construction it takes values in $B_{k+1}$. Let $C_{k}$ be the subset of $B_{k+1}$ of cardinality ${R^d} 4^{d+1}$ satisfying
\begin{equation}
\bigcup_{w\in C_{k}} (w+R_{k}) = B_{k {+ 1}}\cap (\ZZ^d\times L_k \ZZ).
\end{equation}
Our $R$--allowed path $\sigma$ can be divided into four disjoint $R$--allowed paths $\sigma_1,\sigma_2,\sigma_3,\sigma_4$ of length $L_k-1$, starting respectively in $w_1+R_{k}$, $w_2+R_{k}$, $w_3+R_{k}$ and $w_4+R_{k}$, where $w_1,w_2,w_3,w_4\in C_{k}$.

Let us assume by contradiction that we cannot find three among these four allowed paths $\sigma_i$ satisfying
$M(\sigma_i,H)< k^2$. This implies that we can find at least two $i$'s such that $M(\sigma_i,H)\geqslant k^2$. Therefore,
$$M(\sigma,H)=\sum_{i=1}^4 M(\sigma_i,H)\geqslant 2k^2\geqslant (k+1)^2,$$
since $k\geqslant 3$. This contradicts the assumption we made on $\sigma$.

Therefore, for three $i$'s the event $A_{k,H}(w_i)$ occurs. Now, $A_{k,H}(w_i)$ is measurable with respect to $\eta$ restricted to the box $w_i+B_{k}$. In particular two events $A_{k,H}(w_i)$ supported by boxes that are $L_k$-separated in time occur. Also note that by choice of the scale ${L_k}=4^k$ there are at most ${R^d} 4^{d+1}$ possible choices for each of these two boxes (since $w_i\in C_k$). Noting that the invariance of $\eta$ implies that $\bP(A_{k,H}(w))$ does not depend on $w$, the previous argument and the decorrelation hypothesis lead to
\nc{c:decoupling2}
\begin{align}\label{e:inequality_renormalization}
    \P(A_{k+1,H})
    &\leqslant { R^{2d}} 16^{d+1} (\P(A_{k,H})^2+\uc{c:decoupling} L_k^{-\alpha})\nonumber\\
    &= { R^{2d}} 16^{d+1} \P(A_{k,H})^2+\uc{c:decoupling2} L_k^{-\alpha},
\end{align}
where $\uc{c:decoupling2}>0$.

\subsubsection{Triggering}\label{sss:triggering}
In this section we see how to tune the parameter $H$ so that \eqref{e:estimate_A} is satisfied for some fixed scale $L_k$.

Let $k\in\NN$ and $\beta<\alpha$. Let us set $H_k=\lfloor L_k/k^2\rfloor$ and define $C_k'$ to be a minimal subset of $B_k$ satisfying
$$\bigcup_{w\in C_k'} (w+\llbracket 0,H_k-1\rrbracket^d\times\{0\})=B_k\cap (\ZZ^d\times H_k\ZZ).$$
The cardinality of $C_k'$ satisfies $\vert C_k'\vert \leqslant c\,k^{2(d+1)}.$ Then, by Corollary~\ref{cor:threatsinbox},
\begin{align*}
  \P(A_{k,H_k})
  &\leqslant \P\left(\exists w\in C_k',\,\exists z\in w+\llbracket 0,H_k-1\rrbracket^d\times\{0\},\,z\text{ is not $H_k$--threatened}\right)\\
  &\leqslant  ck^{2(d+1)}\uc{c:threatsinbox} H_k^{-\alpha}\\
  &\leqslant c\uc{c:threatsinbox} k^{2(d+1-\alpha)}L_k^{-\alpha}.
\end{align*}

\nc{c:triggering}

At the end of the day, there exists $\uc{c:triggering}=\uc{c:triggering}(\alpha,\beta,\uc{c:decoupling}, \uc{c:eta_zero})>0$ such that
\begin{align}\label{e:triggering}
    \P(A_{k,H_k})\leqslant \uc{c:triggering} L_k^{-\beta}.
\end{align}

\subsubsection{Induction}
\label{sss:induction}

Here we check that, if \eqref{e:estimate_A} holds for a fixed large enough $k$, the estimate propagates to higher values of $k$, thus concluding the proof of Lemma~\ref{l:estimate_A}.

We fix $k_0=k_0({R}, \alpha,\beta,\uc{c:decoupling}, \uc{c:eta_zero})\in\NN$ satisfying
\begin{align}\label{e:condition_k_0}
  { R^{2d}} 16^{d+1} \uc{c:triggering} 4^{-\beta(k_0-1)}+4^{\beta}
  \uc{c:decoupling2}\uc{c:triggering}^{-1} 4^{-(\alpha-\beta) k_0}\leqslant 1.
\end{align}

We fix $H=H_{k_0}=\lfloor L_{k_0}/k_0^2\rfloor$. By \eqref{e:triggering} we have
$$\P(A_{k_0,H})\leqslant \uc{c:triggering} L_{k_0}^{-\beta}.$$
Suppose that $\P(A_{k,H})\leqslant \uc{c:triggering} L_k^{-\beta}$ for some $k\geqslant k_0$. Then, using \eqref{e:inequality_renormalization}, we have
\begin{align*}
  \frac{\P(A_{k+1,H})}{\uc{c:triggering}L_{k+1}^{-\beta}}
  &\leqslant { R^{2d}} 16^{d+1} \P(A_{k,H})^2 \uc{c:triggering}^{-1}
    L_{k+1}^{\beta}+\uc{c:decoupling2}
    L_k^{-\alpha} \uc{c:triggering}^{-1} L_{k+1}^{\beta}\\
  &\leqslant { R^{2d}} 16^{d+1} \uc{c:triggering} 4^{-\beta(k-1)}+4^{\beta}
    \uc{c:decoupling2}\uc{c:triggering}^{-1} 4^{-(\alpha-\beta) k}\\
  &\leqslant { R^{2d}} 16^{d+1} \uc{c:triggering} 4^{-\beta(k_0-1)}+4^{\beta}
    \uc{c:decoupling2}\uc{c:triggering}^{-1} 4^{-(\alpha-\beta) k_0}\\
  &\leqslant 1,
\end{align*}
using \eqref{e:condition_k_0}. This concludes the proof of Lemma~\ref{l:estimate_A} via induction.

\subsection{Moments of $T_1$}
We can now combine Lemmas~\ref{l:fall_on_trap} and~\ref{l:estimate_A}, choosing $\Sigma=\{z\in\ZZ^{d+1},\eta_z=1\}$, to obtain bounds on the moments of $T_1$ (recall its definition in \eqref{e:T_k} and the line underneath).

Recall that $\hat{\bP}$ denotes the law of the random walk together with the environment conditioned on $\eta_{0}=1$.
\begin{lemma}\label{l:moments}
    If $\alpha>0$, and $ p \in (0, \alpha)$, then
    \begin{align}\label{e:moments}
    \hat{\bE}[ T_1^p ] < \infty.
	\end{align}
	In particular, $T_1$ is $\hat{\bP}$-almost surely finite.
	\end{lemma}
\begin{proof}
 The first observation is that
\begin{equation}\label{e:Jtok}
L_k\leqslant JH\Longrightarrow A_{k,H}^c\subseteq\{M_J\geqslant k^2\},
\end{equation}
since in that case all paths of length $JH$ starting at $0$ (which enter in the definition of $M_J$) contain a path $\sigma$ of length $L_k-1$ starting in $R_k$ with $M(\sigma,H)\geqslant k^2$.

Set $\beta \in (p,\alpha)$, fix $H$ as chosen in Lemma~\ref{l:estimate_A} and $J \in \NN$ such that $JH\geq L_{k_0}$. Let us choose $k \in \N$ such that $L_k \leqslant JH < L_{k+1}$; in particular $k \geqslant \frac{\ln J+\ln H}{\ln 4}-1\geqslant c\ln J$ for some universal constant $c>0$.  Now, note that conditioned on $\eta_0=1$, $T_1$ is equal to $T^0$ associated to $\Sigma=\{z\in\ZZ^{d+1},\eta_z=1\}$. Therefore, using Lemma \ref{l:fall_on_trap}, we have
\begin{align*}
\hat{\bP} ( T_1>JH ) & \leqslant \frac{1}{\P(\eta_0=1)} \E\left[(1 - \uc{c:ellipticity}^H)^{M_J}\right]&\mbox{by Lemma~\ref{l:fall_on_trap}} \\
&\leqslant \uc{c:eta_zero}^{-1}\left(\E \big[ (1 - \uc{c:ellipticity}^H)^{M_J}\mathbf{1}_{A_{k,H}^c} \big]+\E\big[ (1 - \uc{c:ellipticity}^H)^{M_J}\mathbf{1}_{A_{k,H}} \big] \right)&\mbox{by \eqref{e:eta_zero}}\\
&\leqslant \uc{c:eta_zero}^{-1}(1 - \uc{c:ellipticity}^H)^{k^2}+\uc{c:eta_zero}^{-1}\P(A_{k,H})&\mbox{by \eqref{e:Jtok}}\\
&\leqslant \uc{c:eta_zero}^{-1}(1-\uc{c:ellipticity}^H)^{(c\ln J)^2}+\uc{c:eta_zero}^{-1}\uc{c:triggering} (H/4)^{-\beta} J^{-\beta}&\mbox{by Lemma~\ref{l:estimate_A}}\\
&\leqslant c J^{-\beta}.
\end{align*}
The result follows from the fact that $\beta>p$ and
\begin{equation}
\hat{\bE}[T_1^p] = \int_0^\infty \hat{\bP} \big( T_1\geqslant x^{1/p} \big) \,\mathrm{d}x,
\end{equation}
finishing the proof of the lemma.
\end{proof}

\subsection{Proof of Theorem \ref{t:main_theorem}}\label{ss:final_proof}

With the tools we developed in the previous sections (decomposition of the trajectory into independent pieces and estimates on the length of the pieces), the proof of Theorem~\ref{t:main_theorem} is essentially standard. We detail it here for completeness.

\subsubsection{Law of large numbers}
We now assume $\alpha>1$ and we study the almost sure convergence of $Z_t/t$. For $t\geqslant 0$, set
\begin{equation}
K_t = \min \big\{ k \geqslant 0, T_k \geqslant t \big\}.
\end{equation}
We then write, for $t>0$,
\begin{equation}\label{e:terms}
  \frac{Z_t}{t} = \frac{Y_{K_t}}{K_t}\frac{K_t}{T_{K_t}} \frac{T_{K_t}}{t} + \frac{Z_t - Y_{K_t}}{t}.
\end{equation}
Note that, since $\alpha>1$, Lemma~\ref{l:moments} yields
\begin{equation}
\hat{\bE} [ \vert Z_{T_1} \vert ] \leqslant R \, \hat{\bE} [T_1] < \infty.
\end{equation}
Therefore, by Lemma \ref{l:indep_decomp}, under $\hat{\bP}$, $(Y_{k+1}-Y_k)_{k \geqslant 0}$ is an i.i.d.\ sequence of integrable random variables with distribution $Z_{T_1}$.

Notice also that Lemma \ref{l:indep_decomp} implies that $\hat{\bP}$-almost surely, $T_k \underset{k\to\infty}{\longrightarrow} \infty$. Therefore $K_t \underset{t\to\infty}{\longrightarrow} \infty$ as well. At the end of the day, the strong law of large numbers implies that, $\hat{\bP}$-almost surely,
\begin{equation}
\frac{Y_{K_t}}{K_t} \underset{t\to\infty}{\longrightarrow} \hat{\bE} [Z_{T_1}].
\end{equation}
In particular, since $Z_{t} = (X_{t}, t)$, it follows that, $\hat{\bP}$-almost surely,
\begin{equation}
\frac{T_{K_t}}{K_t} \underset{t\to\infty}{\longrightarrow} \hat{\bE} [T_1].
\end{equation}

Note now that $\frac{T_{K_t}}{t}=1+\frac{T_{K_t}-t}{t}$, where $\frac{T_{K_t}-t}{t} \leqslant \frac{\vert Z_t-Y_{K_t}\vert}{t}$. Therefore, the only thing left to do is study the convergence of the last term in \eqref{e:terms}, that is, $\frac{Z_t-Y_{K_t}}{t}$. Now, if we fix $\varepsilon>0$ and $\beta = \frac{1+\alpha}{2}$, Markov's inequality yields {(notice that $T_{K_t}-t \geq T_{K_t}-T_{K_t-1}$ for the second inequality)}
\begin{align*}
 \hat{\bP} \left(\frac{\vert Z_t-Y_{K_t}\vert}{t}\geqslant \varepsilon \right)
& \leqslant \hat{\bP} (T_{K_t}-t\geqslant \varepsilon t/R ) \\
& \leqslant \hat{\bP} (T_1 \geqslant \varepsilon t/R ) \leqslant c \, \hat{\bE} \big[ T_1^{\beta} \big] t^{-\beta}.
\end{align*}which is summable in $t$, since $\beta \in (1, \alpha)$ and $\hat{\bE} \big[ T_1^{\beta} \big] < \infty$, by Lemma~\ref{l:moments}. Therefore, using Borel-Cantelli's lemma, $\hat{\bP}$-almost surely,
\begin{equation}
\frac{Z_t-Y_{K_t}}{t} \underset{t\to\infty}{\longrightarrow} 0.
\end{equation}
Putting all the convergences together in~\eqref{e:terms}, we obtain that $\hat{\bP}$-almost surely,
\begin{equation}
\frac{Z_t}{t} \underset{t\to\infty}{\longrightarrow} \frac{\hat{\bE}[Z_{T_1}]}{\hat{\bE}[ T_1 ]}.
\end{equation}

From now on, we let
\begin{equation}\label{d:a}
v = \frac{\hat{\bE}[Z_{T_1}]}{\hat{\bE}[ T_1 ]} \in \RR^{d+1}.
\end{equation}
It remains to show that $\bP$-almost surely, $\frac{Z_t}{t} \underset{t\to\infty}{\longrightarrow} v$. First note that if $T_0<\infty$, for any $t>T_0$, we have
$$\frac{Z_t}{t}=\frac{Z_t-Z_{T_0}}{t-T_0}\,\frac{t-T_0}{t}+\frac{Z_{T_0}}{t},$$
where $\bP$-almost surely, $\mathbf{1}_{T_0<\infty}\frac{t-T_0}{t}\underset{t\to\infty}{\longrightarrow} 1$ and $\mathbf{1}_{T_0<\infty}\frac{Z_{T_0}}{t}\underset{t\to\infty}{\longrightarrow} 0$.

Now, $\mathbf{1}_{T_0<\infty} (Z_{T_0+t}-Z_{T_0})_{t\geqslant 0}$ under $\bP$ has the same distribution as $(Z_t)_{t\geqslant 0}$ under $\hat{\bP}$. Therefore, $\bP$-almost surely,
\begin{equation}
\mathbf{1}_{T_0<\infty} \frac{Z_t-Z_{T_0}}{t-T_0}\underset{t\to\infty}{\longrightarrow} v.
\end{equation}
Finally, since $T_0 \leqslant T_1 < \infty$ almost surely, we obtain the result.

\subsubsection{Central limit theorem}
For the central limit theorem, we strongly rely on \cite{Janson}, following ideas of \cite{Serfozo}. Mind that Theorem 1.4 in \cite{Janson} is stated for real-valued processes, but it can easily be generalized to higher dimensions; when the proof refers to Anscombe's theorem, we use Theorem 3.1 in \cite{gut2009stopped}, that can easily be stated and proven in dimension $d$.

We assume that $\alpha>2$. All assumptions of Theorem 1.4 in \cite{Janson} (considering multidimensional processes) are satisfied (note that we need to apply the theorem with the probability measure $\hat{\bP}$), namely:
\begin{itemize}
 \item $(Z_t)_{t\in\NN}$ has regenerative increments over $(T_n)_{n\in\NN}$, which means that the $$(Z_t-Z_{T_k},\,0\leqslant t\leqslant T_{k+1}-T_k)_{k\in\NN}$$ are i.i.d. Note that this requirement is slightly stronger than Lemma \ref{l:indep_decomp}, but it is clear that it is also satisfied, with the same arguments as in the proof of Lemma \ref{l:indep_decomp}.

 \item The three quantities $\hat{\bE}[T_1]$, $\hat{\bE}[Z_{T_1}]$ and $\mathrm{Cov}(Z_{T_1}- v T_1 \vert \eta_0=1)$ (where $v$ is defined in \eqref{d:a}) are finite, since $\alpha>2$ (recall Lemma \ref{l:moments}).

 \item $\hat{\bP} \Big( \displaystyle\sup_{0 \leqslant t \leqslant T_1} \vert Z_t\vert < \infty \Big) = 1$, which is clear since $\vert Z_t\vert \leqslant Rt$ and $T_1<\infty$ almost surely.
\end{itemize}

Therefore, we have the following convergence in distribution, under $\hat{\bP}$:
\begin{equation}\label{e:cv_distrib}
\frac{Z_t-tv}{\sqrt{t}} \underset{t\to\infty}{\rightharpoonup} \mathcal{N}(0,\Sigma), 
\end{equation}
with $\Sigma = \displaystyle\frac{\mathrm{Cov}(Z_{T_1}-T_1 v \vert \eta_0=1)}{\hat{\bE}[T_1]}$.
In order to conclude, we use the same strategy as for the law of large numbers, by writing, when $T_0<\infty$ and $t>T_0$,
\begin{equation}
\frac{Z_t-tv}{\sqrt{t}}=\frac{Z_t-Z_{T_0}-(t-T_0)v}{\sqrt{t-T_0}}\sqrt{\frac{t-T_0}{t}}+\frac{Z_{T_0}-T_0 v}{\sqrt{t}}.
\end{equation}
Under $\bP$, we have $\mathbf{1}_{T_0<\infty} \frac{Z_t-Z_{T_0}-(t-T_0)v}{\sqrt{t-T_0}} {\rightharpoonup} \mathcal{N}(0,\Sigma)$, $\mathbf{1}_{T_0<\infty} \sqrt{\frac{t-T_0}{t}} \overset{a.s.}{\longrightarrow} 1$, and $\mathbf{1}_{T_0<\infty} \frac{Z_{T_0}-T_0 v}{\sqrt{t}} \overset{a.s}{\longrightarrow} 0$. It follows from Slutsky's theorem that the convergence in \eqref{e:cv_distrib} also holds under $\bP$, which proves the result.

\section{Examples}\label{s:examples}

Let us now collect examples of environments where our main result can be applied. In each case, this amounts to constructing an appropriate field $\eta$, and verifying that it satisfies the Random Markov Property~\eqref{e:rmp} and the decoupling condition from Definition~\ref{d:decoupling}.

\subsection{Boolean percolation}
\nc{c:tail_decay_Boolean}
  Start with a Poisson point process $\Phi$ in $\R^{d+1} \times (0,+\infty)$ with intensity $\Leb(\d x) \otimes \nu(\d \rho)$, where $\nu$ is a fixed distribution in $(0,+\infty)$ that satisfies
  \begin{equation}
    \nu[\rho, \infty) \leqslant \uc{c:tail_decay_Boolean} \rho^{-\beta},
  \end{equation}
for some $\beta>0$, and define the occupied set as the collection of open balls (in the Euclidean norm)
  \begin{equation}
    \mathcal{O} = \bigcup_{(x,\rho) \in \Phi} B(x, \rho) \subseteq\R^{d+1}.
  \end{equation}
  The environment $\omega = (\omega_{z})_{z \in \Z^{d+1}}$ is now given as
  \begin{equation}\label{e:omega_Boolean}
    \omega_{z} = \begin{cases}
      1, & \quad \text{if } z \in \mathcal{O}, \\
      0, & \quad \text{if } z \notin \mathcal{O}.
    \end{cases}
  \end{equation}
  
  As we check below, our results apply to this environment, depending on the exponent $\beta$ that tunes the decay of correlations in this model.
  
  \begin{application}\label{app:boolean}
   If $X$ is a finite-range and uniformly elliptic random walk in environment $\omega$, the conclusions of Theorem~\ref{t:main_theorem} apply for $\alpha=\beta-2d-2$.
  \end{application}

  In order to prove this, let us first define the field $\eta=(\eta_{z})_{z \in \Z^{d+1}}$, consider, for $z \in \Z^{d+1}$, the event
  \begin{equation}
    \mathcal{A}_{z} = \Big\{
    \text{No ball } B(x,\rho) \text{ with } (x,\rho) \in \Phi \text{ intersects both } C_z^{+} \text{ and } C_z^{-}
    \Big\},
  \end{equation}
and set $\eta_{z} = \textbf{1}_{\mathcal{A}_{z}}$. Notice that $\eta_{z}=1$ immediately implies $\omega_{z}=0$, which verifies~\eqref{e:condition_eta}. Furthermore, the collection $(\omega_{z}, \eta_{z})$ inherits the shift invariance property from Boolean percolation.

  Let us now verify that this construction yields a field that satisfies the Random Markov Property \eqref{e:rmp}.
  \begin{lemma}
    The field $(\omega, \eta) = \big( (\omega, \eta)_{z} \big)_{z \in \Z^{d+1}}$ defined as above satisfies the Random Markov Property~\eqref{e:rmp}.
  \end{lemma}

  \begin{proof}
  Let $z\in\ZZ^{d+1}$ and let $W$ be an $\mathcal{F}_z^+$-measurable random variable. We need to prove that, for every $E \in \mathcal{F}_{z}^{-}$,
  \begin{equation}
    \E \big[ W \eta_{z} \textbf{1}_{E} \big] = \E \big[ \E[W|\eta_{z}] \eta_{z} \textbf{1}_{E} \big].
  \end{equation}

  Define now the disjoint sets
  \begin{equation}
    \begin{split}
      A^{+} = \{ (x,t) \in \R^{d+1} \times [0, +\infty): B(x,t) \cap C_{z}^{+} \neq \emptyset \text{ and } B(x,t) \cap C_{z}^{-} = \emptyset \}, \\
      A^{-} = \{ (x,t) \in \R^{d+1} \times [0, +\infty): B(x,t) \cap C_{z}^{-} \neq \emptyset \text{ and } B(x,t) \cap C_{z}^{+} = \emptyset \}, \\
      A^{\pm} = \{ (x,t) \in \R^{d+1} \times [0, +\infty): B(x,t) \cap C_{z}^{+} \neq \emptyset \text{ and } B(x,t) \cap C_{z}^{-} \neq \emptyset \}, \\
    \end{split}
  \end{equation}
  and notice that $W$ is determined by $\Phi \cap (A^{+} \cup A^{\pm})$, $E$ is determined by $\Phi \cap (A^{-} \cup A^{\pm})$ and $\eta_{z} = \textbf{1}_{\Phi \cap A^{\pm} = \emptyset}$. Furthermore, since the sets above are disjoint, $\Phi \cap A^{+}$, $\Phi \cap A^{-}$, and $\Phi \cap A^{\pm}$ are all independent. Note also that on the event that $\Phi\cap A^{\pm}=\emptyset$ (that is, $\eta_z=1$), $W$ can be replaced by $f(\Phi\cap A^+)$ where $f$ is a measurable function. This then implies
  \begin{equation}
    \begin{split}
      \E \Big[ W \eta_{z} \textbf{1}_{E} \Big] & = \E \Big[ W\textbf{1}_{\Phi \cap A^{\pm} = \emptyset} \textbf{1}_{E} \Big] \\
      & = \E \Big[ f(\Phi \cap A^{+}) \textbf{1}_{\Phi \cap A^{\pm} = \emptyset} \textbf{1}_{E} \Big] \\
      & = \E \Big[ \E\ \big[ f (\Phi \cap A^{+}) \big] \textbf{1}_{\Phi \cap A^{\pm} = \emptyset} \textbf{1}_{E} \Big] \\
      & = \E \Big[ \E\ \big[ f (\Phi \cap A^{+})\, \vert\, \eta_z=1\big] \textbf{1}_{\Phi \cap A^{\pm} = \emptyset} \textbf{1}_{E} \Big] \\
      & = \E \Big[ \E\ \big[ W \vert \eta_{z}=1\big] \eta_z \textbf{1}_{E} \Big] \\
      & = \E \Big[ \E\ \big[ W \big| \eta_{z} \big] \eta_z \textbf{1}_{E}\Big].
    \end{split}
  \end{equation}
  In the third equality, we used that $\Phi\cap A^+$ is independent of $\Phi\cap (A^-\cup A^\pm)$. In the fourth equality, we used that $\Phi\cap A^+$ is independent of $\eta$. This concludes the proof.
  \end{proof}

\nc{c:decoupling_boolean_constant}
  To prove Application~\ref{app:boolean}, it only remains to verify that the field $\eta$ satisfies the decoupling property with the correct exponent $\alpha$.
  \begin{proposition}\label{prop:decoupling_boolean}
    There exists $\uc{c:decoupling_boolean_constant} = \uc{c:decoupling_boolean_constant}(\beta, d) >0$ such that the function defined by
    \begin{equation}
      \varepsilon(r,h,s) = \uc{c:decoupling_boolean_constant} r^{d}hs^{-\beta+d+1}
    \end{equation}
  is a decoupling function for $\eta$ defined as above (recall Remark \ref{r:decoupling_function}). In particular, the $\alpha$--decoupling property holds for $\alpha=\beta-2d-2$.
  \end{proposition}

\nc{c:large_difference}
  In order to prove the proposition above, we need a preliminary lemma which requires some additional notation.
  Given $s>0$, define
  \begin{equation}
    \mathcal{A}_{z}^{s} = \bigg\{
      \begin{array}{c}
        \text{No disk } B(x,\rho) \text{ with } (x,\rho) \in \Phi \text{ and } \d(x, z) < \tfrac{s}{2} \\
        \text{ intersects both } C_z^{+} \text{ and } C_z^{-}
      \end{array}
    \bigg\},
  \end{equation}
  and set $\eta^{s}_{z} = \textbf{1}_{\mathcal{A}_{z}^{s}}$. Notice that $\eta^{s}_{z}$ is completely determined by the $\Phi$ restricted to the set $\big( z+ \big(- \tfrac{s}{2}, \tfrac{s}{2} \big)^{d+1} \big) \times [0,+\infty)$.
  We now prove the following lemma
  \begin{lemma}\label{lemma:approximation_boolean}
    Assume $\beta > d$. There exists a constant $\uc{c:large_difference} = \uc{c:large_difference}(\beta, d) >0$ such that, for any $s>0$ and $z \in \Z^{d+1}$,
    \begin{equation}
    \P \big( \eta_{z} \neq \eta_{z}^{s} \big) \leqslant \uc{c:large_difference} s^{-\beta+d+1}.
    \end{equation}
  \end{lemma}

  \begin{proof}
    Assume without loss of generality that $z=0$ and notice that
    \begin{equation}
      \P \big( \eta_{0} \neq \eta_{0}^{s} \big) \leqslant \P \big( \mathcal{B}_{s} \big),
    \end{equation}
    where
    \begin{equation}
      \mathcal{B}_{s} =  \bigg\{
        \begin{array}{c}
          \text{There exists a disk } B(x,\rho) \text{ with } (x,{\rho}) \in \Phi \text{ and } |x| > \tfrac{s}{2} \\
          \text{that intersects both } C_0^{+} \text{ and } C_0^{-}
        \end{array}
      \bigg\}.
    \end{equation}

    In order to bound the probability of $\mathcal{B}_{s}$, we further introduce the event
    \begin{equation}
      \tilde{\mathcal{B}}_{s} =  \bigg\{
        \begin{array}{c}
          \text{There exists a disk } B(x,\rho) \text{ with } (x,\rho) \in \Phi,
          \\ x \in [0,+\infty)^{d+1} \text{ and } |x| > \tfrac{s}{2} \text{ that intersects } C_0^{-}
        \end{array}
      \bigg\},
    \end{equation}
    and notice that rotation invariance implies
    \begin{equation}
      \P \big( \mathcal{B}_{s} \big) \leqslant 2^{d+1} \P \big( \tilde{\mathcal{B}}_{s} \big).
    \end{equation}

\nc{c:lower_bound_distance_cone}
    We now notice that the distance between a point $x \in  [0,+\infty)^{d+1}$ and $C_{0}^{-}$ is lower bounded by the distance from $x$ to the set $\{ (x,t) \in [0, +\infty)^{d} \times (-\infty, 0] : |x| = -Rt \}$, which yields the bound
    \begin{equation}
      \d(x, C_{0}^{-}) \geqslant \uc{c:lower_bound_distance_cone} |x|,
    \end{equation}
    for some $\uc{c:lower_bound_distance_cone}>0$. From this is follows immediately that
    \begin{equation}
      \begin{split}
        \P \big( \tilde{\mathcal{B}}_{s} \big) & \leqslant \int_{[0,+\infty)^{d+1} \setminus B\big( 0,\tfrac{s}{2} \big)} \nu \big[\d(x, C_{0}^{-}), +\infty) \d \, x \\
        & \leqslant \int_{[0,+\infty)^{d+1} \setminus B\big( 0, \tfrac{s}{2} \big)} \nu \big[\uc{c:lower_bound_distance_cone} |x|, +\infty) \d \, x \\
        & \leqslant \int_{[0,+\infty)^{d+1} \setminus B \big( 0, \tfrac{s}{2} \big)} \uc{c:tail_decay_Boolean} \uc{c:lower_bound_distance_cone}^{-\beta} |x|^{-\beta} \d \, x \\
        & \leqslant \uc{c:tail_decay_Boolean} \uc{c:lower_bound_distance_cone}^{-\beta} c \int_{c'\tfrac{s}{2}}^{+\infty} u^{-\beta+d} \d \, u \\
        & \leqslant \tilde{c} s^{-\beta+d+1},
      \end{split}
    \end{equation}
  concluding the proof.
  \end{proof}

  We are now ready to prove Proposition~\ref{prop:decoupling_boolean}.
  \begin{proof}[Proof of Proposition~\ref{prop:decoupling_boolean}]

  Let $A=\prod_{i=1}^{d+1}\llbracket a_i,b_i\rrbracket$ and $B=\prod_{i=1}^{d+1}\llbracket a'_i,b'_i\rrbracket$ be two boxes satisfying \eqref{e:distance_decoupling} and $f_1:\{0,1\}^A\rightarrow\{0,1\}$, $f_2:\{0,1\}^B\rightarrow\{0,1\}$. Consider now the events
  \begin{equation}
    \begin{split}
      G_{A} = \Big\{ \eta_{z} = \eta_{z}^{s}, \text{ for all } z \in A \Big\}, \\
      G_{B} = \Big\{ \eta_{z} = \eta_{z}^{s}, \text{ for all } z \in B \Big\},
    \end{split}
  \end{equation}
and notice that $G_{A}$ and $G_{B}$ are determined by $\Phi$ restricted to the sets $ \prod_{i=1}^{d+1}\big( a_i-\tfrac{s}{2} , b_i +\tfrac{s}{2} \big) \times \big[ 0, +\infty \big)$ and $ \prod_{i=1}^{d+1}\big( a'_i-\tfrac{s}{2} , b'_i +\tfrac{s}{2} \big) \times \big[ 0, +\infty \big)$, respectively. Due to the assumption that $\sep(A,B) \geqslant s$, these two sets are disjoint and thus the configurations $\eta_{z}^{s}$ in the boxes $A$ and $B$ are independent. This then yields
  \begin{equation}
    \begin{split}
    \E[f_1(\eta_A)f_2(\eta_B)] & \leqslant \E \Big[ f_1 \big( \eta_A^{s} \big) f_2(\eta_B^{s}) \Big] + \P(G_{A}) + \P(G_{B}) \\
    & = \E \Big[ f_1 \big( \eta_A^{s} \big) \Big] \E \Big[ f_2 \big( \eta_B^{s} \big) \Big] + \P (G_{A}) + \P( G_{B}) \\
    & \leqslant \E[f_1(\eta_A)]\E[f_2(\eta_B)] + 2\P(G_{A}) + 2\P(G_{B}).
    \end{split}
  \end{equation}

  To conclude, observe that union bound together with Lemma~\ref{lemma:approximation_boolean} implies
  \begin{equation}
    \P(G_{A}) \leqslant \uc{c:large_difference} r^{d}h s^{-\beta+d+1}.
  \end{equation}
  An analogous bound is satisfied by $\P(G_{B})$, concluding the proof.
  \end{proof}

\subsection{Independent renewal chains}

We now assume that the environment $\omega$ is composed of independent copies of a renewal chain than we introduce now. Fix an aperiodic\footnote{We say that $\mu$ is aperiodic if $\gcd\{k \in \Z_{+}: \mu(k+1) >0\} = 1$.} reference distribution $\mu$ on the non-negative integers with finite expectation. For each $x \in \Z^{d}$, the environment $\omega(x, \cdot\,)$ will be given by an independent copy of a renewal process, defined as the Markov chain with transition probabilities $p(k,k-1) =1$, if $k \geqslant 1$, and $p(0,k) = \mu(k)$. We call $\mu$ the \emph{interarrival distribution} of the renewal chain and assume that the initial state $\omega(x,0)$ is such that each of these chains is stationary. This amounts to taking $\omega(x,0)$ with distribution $\hat\mu$ defined as
\begin{equation}\label{eq:initial_distribution_renewal}
\hat\mu( k)= \frac{1}{\E[\xi]+1} \P(\xi \geqslant k), \text{ with }\xi\sim\mu.
\end{equation}

Let us note for future reference that, if $\xi \sim \mu, \hat{\xi} \sim \hat{\mu}$ and $\E(\xi^{1+\beta}) < \infty$ for some $\beta > 0$, then
\begin{equation}\label{eq:moment}
\E[ \hat{\xi}^{\beta}] < \infty.
\end{equation}

Let us also recall the following theorem from~\cite{zbMATH03604142}.
\begin{theorem}[\cite{zbMATH03604142}]\label{t:coupling_renewal}
Let $\P_{\mu}^{\nu, \nu'}$ denote the distribution of two independent renewal chains ${ V}$ and ${ \tilde{V}}$ with independent initial distributions $\nu$ and $\nu'$ and same interarrival distribution $\mu$. Denote the coupling time of these chains by
\begin{equation}
T = \inf\{ n \geqslant 0: { V_{n}} = { \tilde{V}_{n}} = 0 \}.
\end{equation}
\begin{enumerate}
 \item If $\mu$ is aperiodic and has finite $1+\beta$ moment for some $\beta>0$, and both $\nu$ and $\nu'$ have finite $\beta$ moment, then $\E_{\mu}^{\nu, \nu'} \big[ T^{\beta} \big] < \infty$.
 \item If $\mu$ is aperiodic and has finite $1+\beta$ moment for some $\beta\geqslant 0$, and both $\nu$ and $\nu'$ have finite $1+\beta$ moment, then $\E_{\mu}^{\nu, \nu'} \big[ T^{1+\beta} \big] < \infty$.
\end{enumerate}
\end{theorem}

Our control of the correlations in the environment will be based on this theorem. In particular, it depends mostly on the exponents $\beta$ for which \eqref{eq:moment} holds. 

\begin{application}
 Assume\footnote{There are non-trivial distributions for which $\gamma_{\mu} = 0$ (take $\mu(0)=0$ or $\mu(2k)\propto 2^{-k}$, $\mu(2k+1)\propto k^{-2}$). On the other hand, notice that any non-trivial $\mu$ for which $\gamma_{\mu} > 0$ is necessarily aperiodic.} that $\mu$ is such that
 \begin{equation}\label{eq:gamma}
  \gamma_{\mu} := \inf_{k\colon \hat\mu(k)> 0}\frac{\mu(k)}{\hat\mu(k)}>0.
\end{equation}
If $X$ is a finite-range and uniformly elliptic random walk in environment $\omega$, and $\mu$ has a finite moment of order $1+\beta>1+2+2d$  the conclusions of Theorem~\ref{t:main_theorem} apply for $\alpha=\beta-2-2d$.
\end{application}

We now propose an explicit construction of the environment $\omega$ that will allow us to construct an appropriate field $\eta$ encoding the fact that $\omega$ sometimes ``forgets its past'', by looking for times at which it couples well with an independent stationary version of itself.

\subsubsection{Graphical construction}\label{s:graphical-construction-renewal}

Assume $\gamma_{\mu} > 0$. Consider three i.i.d.\ independent collections $(\hat W_t^x)_{x\in\Z^d,t\in\Z}$, with distribution $\hat\mu$, $(Z^x_t)_{x\in\Z^d,t\in\Z}$, with distribution $\mathrm{Ber}(\gamma_{\mu})$, and $(Y^x_t)_{x\in\Z^d,t\in\Z}$, with distribution\footnote{$\gamma_{\mu} < 1$ unless $ \mu=\delta_0$; in the latter case the distribution of $Y$ is irrelevant and can be chosen arbitrarily.} $\frac{1}{1-\gamma_{\mu}}(\mu-\gamma_{\mu} \hat{\mu})$. Note that the latter is indeed a probability measure on the non-negative integers, thanks to the definition of $\gamma_{\mu}$. We now define
\begin{equation}
W_t^x:=Z_t^x\hat W^x_t+(1-Z_t^x)Y_t^x.
\end{equation}
The $W_t^x$ are i.i.d.\ with distribution $\mu$. Also note, that if $Z_t^x=1$, $W_t^x$ is equal to a random variable with stationary distribution $\hat\mu$ independent from the past of the chain at $x$.

We aim to construct $(\omega(x,t))_{x\in\Z^d,t\in\Z}$ such that
\begin{itemize}
 \item the processes $(\omega(x,\cdot))_{x\in\Z^d}$ are independent;
 \item for all $x\in\Z^d$, $(\omega(x,t))_{t\in\Z}$ is a stationary renewal chain that uses the collection $W^x$ to determine its jumps, i.e. $\omega(x,t)=W_t^x$ if $\omega(x,t-1)=0$.
\end{itemize}
By independence, it suffices to fix $x\in\Z^d$ and construct $(\omega(x,t))_{t\in\Z}$ satisfying the second property. For $K\in\N$, define ${ V}^{(K)}$ as the renewal chain started from $0$ at time $-K$ whose jumps are determined by $W^x$: for $t \geqslant -K$,
\begin{equation}
\begin{split}
  { V}^{(K)}_{-K}
  & = 0, \\
  { V}^{(K)}_{t+1}
  & =
    \begin{cases}
      { V}^{(K)}_t-1 & \text{if } { V}^{(K)}_t \geqslant 1, \\
      W^x_{t+1} & \text{if } { V}^{(K)}_t = 0.
    \end{cases}
\end{split}
\end{equation}

\begin{lemma}
 If $\mu$ is aperiodic and has finite moment of order $1+\delta$, for some $\delta>0$, 
 \begin{enumerate}
  \item for any $t\in\Z$, $\omega(x,t) := \lim_{K\rightarrow \infty} { V}^{(K)}_t $ exists a.s.;
  \item the process $\omega(x,\cdot)$ thus obtained is a stationary renewal chain whose jumps are determined by $W^x$. 
 \end{enumerate}
\end{lemma}
\begin{proof}
 Since our construction is invariant under time translations, it is enough to prove the first item for $t=0$. We claim that there exists $c\in\R_+$ such that
 \begin{equation}
  \P({ V}_0^{(K)}\neq { V}_0^{(K+1)})\leqslant cK^{-(1+\delta)}.
 \end{equation}
Indeed, let $T:=\inf\{k\geqslant 0\colon { V}_{-K+k}^{(K)} = { V}_{-K+k}^{(K+1)}\}$ denote the coupling time of the chains started at time $-K$ from $0$ and from ${ V}^{(K+1)}_{-K}=W^x_{-K}$. Note that, up to time $T$, the chains $({ V}^{(K)}_{-K+t})_{t\geqslant 0}$ and $({ V}^{(K+1)}_{-K+t})_{t\geqslant 0}$ are independent (they are constructed using disjoint sets of independent variables). By Theorem~\ref{t:coupling_renewal} {part 2.}, $T$ has finite $1+\delta$ moment. Also note that the distribution of $T$ does not depend on $K$.

It remains to notice that
\begin{equation}
\P\big( { V}_0^{(K)} \neq { V}_0^{(K+1)} \big) \leqslant \P(T>K)
\end{equation}
and to apply Markov's inequality.

To show that $({ V}_0^{(K)})_K$ converges almost surely, let us now estimate
\begin{equation}
\begin{split}
\P(\{\exists K_0\text{ s.t. }\forall K\geqslant K_0,\ { V}_0^{(K)}={ V}_0^{(K_0)}\}^c)&\leqslant \P(\forall K_0,\ \exists K\geqslant K_0\text{ s.t. }{ V}_0^{(K)}\neq { V}_0^{(K+1)})\nonumber\\
&\leqslant \inf_{K_0}\sum_{K\geqslant K_0}\P({ V}_0^{(K)}\neq { V}_0^{(K+1)})\\
&\leqslant c\inf_{K_0}\sum_{K\geqslant K_0}\frac{1}{K^{1+\delta}}=0,
\end{split}
\end{equation}
since $\delta>0$. The first item is proved.

It is now immediate to check that $\omega(x,\cdot)$ is a renewal chain whose jumps are determined by $W^x$: by our construction, if $({ V}_{t_0}^{(K)})_{K\geqslant K_0}$ is constant, so is $({ V}_{t}^{(K)})_{K\geqslant K_0}$ for any $t\geqslant t_0$. Therefore, for all $t_0\in\Z$ there exists a random $K$ s.t. $(\omega(x,t))_{t\geqslant t_0}=({ V}_t^{(K)})_{t\geqslant t_0}$.

To show that $\omega(x,0)\sim\hat\mu$, consider $\hat { V}^{(K)}$ the (stationary) chain started at time $-K$ from an independent variable $\hat { V}^{(K)}_{-K}\sim\hat\mu$ and using the $W^x$ variables as jumps. Then
\begin{equation}
\begin{split}
 \P \big( \omega(x,0) \neq \hat { V}_0^{(K)} \big) & \leqslant \P\big( \omega(x,0) \neq { V}_0^{(K)} \big) + \P \big( { V}_0^{(K)} \neq \hat { V}_0^{(K)} \big) \\
 & \leqslant \P\big( \omega(x,0) \neq { V}_0^{(K)} \big)+\P( \hat{T} > K),
\end{split}
\end{equation}
where $\hat{T}$ denotes the coupling time between $({ V}^{(K)}_t)_{t \geqslant -K}$ and $(\hat{{ V}}^{(K)}_t)_{t \geqslant -K}$. By Theorem~\ref{t:coupling_renewal} {part 1.}, $\hat{T}$ has finite moment of order $\delta$, implying that the right-hand side above tends to $0$ as $K$ grows.
\end{proof}

\subsubsection{Construction and properties of $\eta$}

We now build $\eta$. For this, we look for space-time points so that, somewhere between the cones of the past and future trajectories rooted at $x$, $\omega$ couples with a version of itself that is independent from its past.

Let us fix $(x_0,t_0)\in\Z^d\times \Z$ and construct $\eta(x_0,t_0)$. Let $t^x := -\lceil |x|/R \rceil +1$ and consider the random times (recall the definition of the $Z_t^x$ from the previous section)
\begin{equation}\label{e:defT_x}
  T^x = T^x(x_0,t_0) := \inf\{t \geqslant t_0+t^{(x {- x_0})} : \omega(x,t)=0 \text{ and } Z_{t+1}^x = 1 \}.
\end{equation}

\begin{definition}
 Let $\eta(x_0,t_0)$ be the indicator function of the event
\begin{equation}\label{e:defeta-renewal}
 \{\forall x\in\Z^d\setminus\{x_0\},\ T^x<t_0+|x-x_0|/R,\text{ and }T^{x_0}=t_0\}.
\end{equation}
\end{definition}

We now verify that, with this construction, $(\omega,\eta)$ satisfies the hypotheses of our theorem. Note that $\eta$ is not exactly measurable w.r.t.\@ $\omega$ itself, but rather w.r.t.\@ the random variables which we use to sample $\omega$. In other words, to apply Theorem~\ref{t:main_theorem}, we need to consider that the environment is given by $(\omega(x,t),\hat W_t^x,Z_t^x,Y_t^x)_{(x,t)\in\Z^d\times\Z}$.

Let us first estimate the probability that $\eta(x_0,t_0)=1$. 

\begin{lemma}
 If $\mu$ is a non-trivial distribution with $\gamma_{\mu} > 0$ (see~\eqref{eq:gamma}) and finite moment of order $1+\beta$, for some $\beta>d$, then 
 \begin{equation}
  \P( \eta(x_0,t_0) = 1 )>0.
 \end{equation}
\end{lemma}
\begin{proof}
 By translation invariance we can assume $x_0=t_0=0$. By independence of the renewal chains at different sites,
 \begin{equation}
 \begin{split}
   \P(\eta(0,0) = 1)
   & =\P(T^0=0)\prod_{x\neq 0}\P(T^x<|x|/R)\\
   & \geqslant \P(T^0=0)\prod_{\substack{|x| \leqslant R/2 \\ x \neq 0}}
   \P\Big(T^x<\frac{|x|}{R}\Big) \prod_{|x|>R/2}
   \P \Big( T^x-t^x < 2\frac{|x|}{R}-1 \Big).\\
   & { \geqslant \P(T^0=0)\prod_{\substack{|x| \leqslant R/2 \\ x \neq 0}}
   \P\Big(T^x<\frac{|x|}{R}\Big)
   \exp \Big\{ - c \sum_{r \geq R/2} r^{d - 1}
   \P \Big( T^x-t^x \geq 2\frac{r}{R}-1 \Big) \Big\}}.
 \end{split}
 \end{equation}
Let us reformulate the distribution of $T^x-t^x$. Let $(Z_i)_{i\in\N}, (Y_i)_{i\in\N}, \hat W$ be independent random variables, with $Z_i \sim \mathrm{Ber}(\gamma_{\mu})$, $Y_i \sim \frac{1}{1-\gamma_{\mu}}(\mu-\gamma_{\mu} \hat{\mu})$, $\hat W \sim \hat{\mu}$. Let $S=\inf\{i \geqslant 1 \colon Z_i=1\}$. Then
\begin{align}\label{e:distT_x}
T^x-t^x\sim \hat W+\sum_{i=1}^{S-1}(Y_i+1).
\end{align}
{ Notice that $W$ and $Y_i$ have finite moment of order $\beta$. Now, Cauchy-Schwartz inequality implies
\begin{equation}
\E \Big[ \Big(\sum_{i=1}^{S-1}(Y_i+1) \Big)^{\beta} \Big] \leq \E \Big[ (S-1)^{\beta-1} \sum_{i=1}^{S-1}(Y_i+1)^{\beta} \Big] = \E[(S-1)^{\beta}] \E[(Y_{1}+1)^{\beta}],
\end{equation}
the last equality following from conditioning on the value of $S$.} In particular, $T^x-t^x$ has finite moment of order $\beta$. 

Since $\beta>d$, by Markov inequality, the last product converges, and it is not difficult to check that every term in the (finite) remaining product is positive, concluding the proof.
\end{proof}

\bigskip

We now need to check that $\eta$ satisfies the decoupling property in Definition~\ref{d:decoupling}. To that end, we proceed as in the case of Boolean percolation, by introducing finite range approximations of $\eta$. For any integer $s \geqslant 1$, let $\eta^{s}(x_0,t_0)$ be the indicator function of the event
\begin{equation}
 \{\forall x \in \Z^d\setminus\{x_0\} \text{ s.t. } |x-x_0| \leqslant s, T^x < t_0 + |x-x_0|/R, \text{ and } T^{x_0}=t_0 \}.
\end{equation}
\nc{c:finite_range_approximation}

\begin{lemma}\label{lemma:finite_range_approximation_renewal}
 Assume that $\mu$ is a non-trivial distribution with $\gamma_{\mu} > 0$ and finite moment of order $1+\beta$, for some $\beta>d+1$. There exists $\uc{c:finite_range_approximation}>0$ such that, for any integer $s \geqslant 1 $ and $(x,t) \in \Z^{d} \times \Z_{+}$,
\begin{equation}
\bP\big( \eta(x,t) \neq \eta^{s}(x,t) \big) \leqslant \uc{c:finite_range_approximation}s^{d-\beta+1}.
\end{equation}
\end{lemma}
\begin{proof}
 By translation invariance we can assume $x=t=0$. By union bound, for $s>R/2$,
 \begin{align}
  \bP\big( \eta{(0,0)} \neq \eta^{s}{(0,0)}  \big) &\leqslant \sum_{|x|>s}{\P(T^x-t^x\geqslant 2|x|/R)}\\
  &\leqslant c\sum_{|x|>s}\frac{1}{|x|^{\beta}}\leqslant c'\sum_{k>s}\frac{1}{k^{\beta-d}}\leqslant c''s^{d-\beta+1},
 \end{align}
where, in the second inequality, we used Markov inequality {with the $\beta$-th power of $T^x-t^x$ and the fact that $T^x-t^x$ has finite $\beta$ moment not depending on $x$ (see \eqref{e:distT_x} and below)}. It then remains to adjust the value of $\uc{c:finite_range_approximation}$ to accommodate $s \leqslant R/2$.
 
\end{proof}

\nc{c:renewal_decoupling}
In order to verify the decoupling property, we {use} an analogous result that holds for one-dimensional renewal chains. The proof is basically the same as that of \cite[Lemma 2.1]{HSS19}.
\begin{lemma}\label{lemma:1d_renewal_decoupling}
Assume $\mu$ is a non-trivial distribution with $\gamma_{\mu} > 0$ and finite moment of order $1+\beta$, for some $\beta>0$. Consider a stationary renewal process $({ V}_{n})_{n \in \Z}$ with interarrival distribution $\mu$. There exists a positive constant $\uc{c:renewal_decoupling}$ such that, for all $m\in\Z,n \in \Z_{+}$ and any pair of events $A$ and $B$ with
\begin{equation}
A \in \sigma( { V}_{k}: k \leqslant m ) \quad \text{and} \quad B \in \sigma( { V}_{k}: k \geqslant m+n),
\end{equation}
it holds that
\begin{equation}
\P(A \cap B) \leqslant \P(A)\P(B) + \uc{c:renewal_decoupling} n^{-\beta}.
\end{equation}
\end{lemma}

We are now in position to prove the decoupling property for $\eta$.
\nc{c:decoupling_renewal_2}
\begin{proposition}
Assume that $\mu$ is a non-trivial distribution with $\gamma_{\mu} > 0$ and has finite moment of order $1+\beta$, for some $\beta>d+1$. 
There exists a positive constant $\uc{c:decoupling_renewal_2}$ such that the function
\begin{equation}
\varepsilon(r,h,s) = \uc{c:decoupling_renewal_2}(r+s)^{d} h s^{d-\beta+1},
\end{equation} 
is a decoupling function for the field $\eta$ (see Remark~\ref{r:decoupling_function}).
\end{proposition}

\begin{proof}
Let $A=\prod_{i=1}^{d+1}\llbracket a_i,b_i\rrbracket$ and $B=\prod_{i=1}^{d+1}\llbracket a'_i,b'_i\rrbracket$ be two boxes satisfying \eqref{e:distance_decoupling} and consider two measurable functions $f_1:  \{0,1\}^{A} \rightarrow \{0,1\}$, $f_2:  \{0,1\}^{B} \rightarrow \{0,1\}$. Assume without loss of generality that $b_{d+1} \leqslant a'_{d+1}$.

Lemma~\ref{lemma:finite_range_approximation_renewal} implies
\begin{equation}\label{eq:prob_difference}
\begin{split}
\P \big( f_{1}(\eta_{|A}) \neq f_{1}(\eta_{|A}^{s/3}) \big) & \leqslant \P \big( \eta_{x,t} \neq \eta_{x,t}^{s/3}, \text{ for some } (x,t) \in A \big) \\
& \leqslant \uc{c:finite_range_approximation}{3}^{-d+\beta-1}r^{d}hs^{d-\beta+1}.
\end{split}
\end{equation}
Observe now that $f_1(\eta_{|A}^{s/3})$ depends only on the fields $\omega$ and $W$ restricted to the set $A' = \prod_{i=1}^{d+1}\llbracket a_i - \frac{s}{3}, b_i+\frac{s}{3} \rrbracket$. Similarly, 
\begin{equation}
 \P \big( f_{2}(\eta_{|B}) \neq f_{2}(\eta_{|B}^{s/3}) \big)  \leqslant\uc{c:finite_range_approximation}{3}^{-d+\beta-1}r^{d}hs^{d-\beta+1},
\end{equation}
and $f_2(\eta_{|B}^{s/3})$ depends only on the fields $\omega$ and $W$ restricted to $B' = \prod_{i=1}^{d+1}\llbracket a'_i - \frac{s}{3}, b'_i+\frac{s}{3} \rrbracket$. In particular, by successive conditioning, Lemma~\ref{lemma:1d_renewal_decoupling} yields
\begin{equation}
 \E\big[f_1(\eta_{|A}^{(s/3)})f_2(\eta_{|B}^{(s/3)})\big]\leqslant \E\big[f_1(\eta_{|A}^{(s/3)})\big]\E\big[f_2(\eta_{|B}^{(s/3)})\big]+\uc{c:renewal_decoupling}3^{\beta} (r+s)^{d} s^{-\beta}.
\end{equation}
This, when combined with~\eqref{eq:prob_difference}, implies
\begin{multline}
 \E[f_1(\eta_{|A})f_2(\eta_{|B})]\leqslant \E[f_1(\eta_{|A})]\E[f_2(\eta_{|B})]\\+\uc{c:renewal_decoupling}3^{\beta} (r+s)^{d} s^{-\beta}+4\uc{c:finite_range_approximation}{3}^{-d+\beta-1}r^{d}hs^{d-\beta+1},
\end{multline}
concluding the proof.
\end{proof}

\bigskip

It remains to prove that the field satisfies the Random Markov Property~\eqref{e:rmp}.
\begin{proposition}
If $\mu$ is a non-trivial distribution with $\gamma_{\mu} > 0$ and finite moment of order $1+\beta$, for some $\beta>d$, the field $(\omega, \eta)$ satisfies the Random Markov Property~\eqref{e:rmp}.
\end{proposition}

\begin{proof}
Fix $z \in \Z^{d}\times \Z_{+}$, a $\mathcal{F}_{z}^{+}$-measurable function $f(\omega_{|C^+_z})$ and an event $E \in \mathcal{F}_{z}^{-}$. Our goal is to prove that
\begin{equation}\label{e:rmpproof}
\E\big[ f(\omega_{|C^+_z}) \eta_{z} \textbf{1}_{E}\big] =  \E\big[ \E[ f(\omega_{|C^+_z}) | \eta_{z}] \eta_{z} \textbf{1}_{E}\big]
\end{equation}
By space-time translation invariance, we may and do assume $z=0$. We also assume that $f$ has finite space support: there exists $\Lambda\subset\Z^d$ such that $f$ depends only on $\omega$ restricted to $C_0^+\cap(\Lambda\times \Z)$. The general claim follows from this particular case through standard arguments. 

Let us introduce some notation. For $x\in\Z^d$, $t\in\Z$, let 
\begin{equation}
\mathcal G_t^x = \sigma \left( \hat{W}_s^x, Z_{s+1}^x, Y_s^x ; s \leqslant t \right),
\end{equation}
and notice that $T^x (= T^x(0,0))$ defined in~\eqref{e:defT_x} is a finite stopping time for the filtration $(\mathcal G_t^x)_t$. Let $\mathcal G_{T^x}^x$ be the associated stopped $\sigma$-algebra and $\mathcal G={\otimes}_{x\in\Z^d} \mathcal G_{T^x}^x$. 
Finally, we define
\begin{equation}
  \mathcal L=\{(x,T^x), x\in\Lambda\},
\end{equation}
on the event $\eta_0=1$, and $\mathcal L=\emptyset$ else. Let $\mathfrak L$ be the set of possible non-empty values for $\mathcal L$.

For $L=\{(x,t_x),x\in\Lambda\}\in\mathfrak L$, we define $\mathcal G_L:=(\otimes_{x\in\Lambda}\mathcal G_{t_x})\otimes(\otimes_{x\in\Lambda^c}\mathcal G_{T^x}^x)$, and $\mathcal H_L=\otimes_{x\in\Lambda}\sigma\left(\hat W_s^x,Z_{s+1}^x,Y_s^x;s> t_x\right)$.

The following claims imply \eqref{e:rmpproof}:
\begin{enumerate}
 \item $\eta_0=\sum_{L\in\mathfrak L}\textbf{1}_{\mathcal L=L}$;
 \item for $L\in\mathfrak L$, $\{\mathcal L=L\}$ is $\mathcal G_L$-measurable;
 \item for $L\in\mathfrak L$, on $\{\mathcal L=L\}$, $\textbf{1}_E$ coincides with a $\mathcal G_L$-measurable variable $ G^L$; 
 \item for $L\in\mathfrak L$, on ${ \{ \mathcal L=L \}}$, $\omega_{|C_0^+\cap(\Lambda\times\Z)}$ coincides with a $\mathcal H_L$-measurable random variable $\omega^{L}$ whose distribution does not depend on $L$.
\end{enumerate}

Let us check that we reach the desired conclusion before proving the claims. 
\begin{align}
 \E[ f(\omega_{|C_0^+}) \eta_0\textbf{1}_E]&=\sum_{L\in\mathfrak L}\E[ f(\omega_{|C_0^+}) \textbf{1}_{\mathcal L=L}\textbf{1}_E]\\
 &=\sum_{L\in\mathfrak L}\E[f(\omega^L)\textbf{1}_{\mathcal L=L}G^L]\\
 &=\sum_{L\in\mathfrak L}\E[ f(\omega_{|C_0^+}) ]\E[\textbf{1}_{\mathcal L=L}\textbf{1}_E]\\
 &=\E[ f(\omega_{|C_0^+}) ]\E[\eta_0\textbf{1}_E],
\end{align}
where the first claim was used in the first and last lines, the other claims in the second line. In the third line, we made use of the independence of $\mathcal G_L$ and $\mathcal H_L$ as well as the third and fourth claims. By taking $E$ as the whole sample space above we get
\begin{equation}
\begin{split}
\E[ f(\omega_{|C_0^+})|\eta_0=1] & = \P(\eta_0=1)^{-1} \E[ f(\omega_{|C_0^+}) \textbf{1}_{\eta_{0}=1}] \\
 & = \P(\eta_0=1)^{-1} \E[ f(\omega_{|C_0^+}) \eta_{0}] \\
 & = \P(\eta_0=1)^{-1} \E[ f(\omega_{|C_0^+})] \E[\eta_{0}] \\
 & = \E[ f(\omega_{|C_0^+}) ],
\end{split}
\end{equation}
so that, since $\E[ f(\omega_{|C_0^+})|\eta_0]\eta_0=\E[ f(\omega_{|C_0^+})|\eta_0=1]\eta_0$,
\begin{equation}
 \E[\E[ f(\omega_{|C_0^+})|\eta_0]\eta_0\textbf{1}_E]=\E[ f(\omega_{|C_0^+}) ]\E[\eta_0\textbf{1}_E].
\end{equation}

Let us now prove the claims. The first and second ones are clear from the definition of $\mathcal L$ and $\mathcal G_L$, let us focus on the third one. Let $L = \{(x,t_x),x\in\Lambda\} \in \mathfrak L$. $\textbf{1}_E$ is a function of $(\omega_z,\eta_z)_{z\in C_0^-}$. By definition of $\mathfrak L$, $\omega_z$ is $\mathcal G_L$-measurable for $z \in C_0^-$. Let us check that $\eta_z$ coincides with a $\mathcal G_L$-measurable random variable for $z=(x_0,t_0)\in C_0^-$. Let $\tilde T^x=T^x(x_0,t_0)$ as in~\eqref{e:defT_x}. Since $\tilde T^x\leqslant T^x$ for any $x\in\Lambda ^c$, $\{\tilde T^x<t_0+|x-x_0|/R\}$ is $\mathcal G_{T^x}^x$-measurable. For $x\in\Lambda$, either $t_x \geqslant t_0+|x-x_0|/R$, in which case $\{\tilde T^x<t_0+|x-x_0|/R\}$ is $\mathcal G_{t_x}$-measurable; or $t_x< t_0+|x-x_0|/R$ and on $\mathcal L=L$ the events $\{\tilde T^x<t_0+|x-x_0|/R\}$ and $\{\tilde T^x<t_x\}$ coincide since $\tilde{T}^x \leqslant T^x$.

Finally, it remains to verify the fourth claim. On $\{\mathcal L=L\}$, $\omega_{|C_0^+}$ coincides with the process started from the $\hat W_{T^x+1}^x$ on $\{(x,T^x+1),x\in\Lambda\}$ which uses $(W^x_t)$ with $(x,t) \in \{(x,t) \colon x \in \Lambda, t> T^x+1\}$ as jumps. This process is stationary and, in particular, its distribution in $C_0^+$ does not depend on $L$. It is also clearly $\mathcal H_L$-measurable.

This concludes the proof of the claims, and thus of the Random Markov Property.
\end{proof}

\section{Open problems}

We finish the paper with two simple examples of dynamics for which the techniques presented have not proved sufficient to establish a CLT.
Namely, we explain in detail Finitary Factors of i.i.d.\ fields and Kinetically Constrained Models.
Other interesting open problems are: the Renewal Contact Process \cite{fontes2023renewalcontactprocessesphase} and Gaussian Fields with fast decaying covariances \cite{baldasso2023fluctuationboundssymmetricrandom}.

\paragraph{Finitary factors of i.i.d.}

A random environment $(\omega_z)_{z \in \mathbb{Z}^{d + 1}}$ is said to be a \emph{factor of i.i.d.} if there exists a measurable function $f: [0, 1]^{\mathbb{Z}^{d + 1}} \to S$ such that $\omega_z = f(V \circ \theta_z), \text{ for all $z \in \mathbb{Z}^{d + 1}$}$ where $V = (V_{z'})_{z' \in \mathbb{Z}^{d + 1}}$ is an i.i.d.\ field of uniform $[0, 1]$ variables and $\theta_z$ stands for the shift by~$z$.

One way to further improve its mixing conditions is to impose that it be a \emph{finitary factor of i.i.d.} (ffiid).
Intuitively speaking, we say $\omega$ is a ffiid if for every $v \in [0, 1]^{\mathbb{Z}^{d + 1}}$, there exists a finite (and random) coding radius $\rho < \infty$, such that knowing the values of the uniform variables inside the ball $B(0, \rho)$, fully determines the value of $f$.
See \cite{10.1214/16-AOP1127} for a very interesting application of this definition on the subject of proper colorings of $\mathbb{Z}^d$.

Assuming that the coding radius $\rho$ has light tails, can we prove a CLT for the random walk on top of a ffiid field? Our attempts to find an appropriate definition of $\eta$ that would satisfy Definition~\ref{d:rmp} have failed.

\paragraph{Kinetically constrained models}

Kinetically constrained models (KCM) are a class of interacting particle systems whose non-equilibrium study poses many challenges. They evolve under (reversible) Glauber dynamics, subject to a constraint that forbids any update unless a certain neighborhood is free of particles. We refer to \cite{KCSM} for a general definition. The most emblematic example of such systems is the East model, which in dimension $1$ can be defined through its generator
\begin{equation}
 \label{e:generatorEast}
 \mathcal L f(\sigma)=\sum_{x\in\Z}(1-\sigma(x+1))(p(1-\sigma(x))+(1-p)\sigma(x))\left[f(\sigma^x)-f(\sigma)\right],
\end{equation}
where $p\in(0,1)$ is a density parameter and $\sigma^x$ denotes the configuration $\sigma\in\{0,1\}^\Z$ flipped at~$x$.

KCM provide a class of environments that display fast, but non-uniform mixing. For instance, the East model is exponentially mixing in $L^2(\mu_p)$ at any density ($\mu_p$ being the product Bernoulli measure with density $p$), but the constraint makes it non-uniformly mixing. They are also generically non-attractive. As a result, few criteria implying limit theorems for RWRE apply to them, although it is expected that at least LLN and CLT should follow from the fast mixing.

In \cite{Avena2017}, law of large numbers and CLT were proved for random walks on KCM with exponential mixing (e.g. the East model at any density), but only in a perturbative regime where the jumps rates of the random walker are close to rates that would make the environment seen from the walker stationary w.r.t.\@ $\mu_p$. The results in \cite{cavalinhos,allasia2023law} show LLN for random walks on the one-dimensional East model assuming only finite-range jumps. However, this is still limited to dimension $1$, and does not yield a CLT.

Can one prove a CLT for random walks on top of the East Model over all range of parameters?

\printbibliography
\end{document}